\documentclass[12pt, reqno, a4paper,pagebackref]{amsart}

\usepackage{amssymb,amsmath,dsfont}
\usepackage{pdfsync}
\usepackage{comment}
\usepackage{stackengine}
\usepackage{graphicx,fancyhdr}
\usepackage{geometry}
\usepackage{color}

\usepackage[colorlinks=true,linkcolor=blue,citecolor=magenta]{hyperref}

\usepackage{cleveref}

\newtheorem{df}{Definition}[section]
\newtheorem{lemma}[df]{Lemma}
\newtheorem{prop}[df]{Proposition}
\newtheorem{thm}[df]{Theorem}

\newtheorem{cor}[df]{Corollary}
\newtheorem{remark}[df]{Remark}

\makeatletter \@addtoreset{equation}{section}

\providecommand{\keywords}[1]{\textbf{\textit{Keywords:}} #1}
\newcommand{\bbP}{{\mathbb{P}}}
\newcommand{\bbE}{{\mathbb E}}
\newcommand{\N}{{\mathbb{N}}}
\newcommand{\R}{{\mathbb{R}}}
\newcommand{\fT}{{\mathfrak T}}
\newcommand{\T}{{\mathbb{T}}}
\newcommand{\eps}{\varepsilon}
\newcommand{\om}{\omega}
\newcommand{\la}{\lambda}
\newcommand{\si}{\sigma}

\newcommand{\cal}{\mathcal}

\makeatother

\begin{document}

\title[Anomalous diffusion limit]
  {Anomalous diffusion limit for a kinetic equation with~a~thermostatted
  interface}

\author{Krzysztof Bogdan}
 \thanks{K.B. was supported through the DFG-NCN Beethoven Classic 3 programme, contract no. 2018/31/G/ST1/02252 (National Science Center, Poland) and SCHI-419/11–1 (DFG, Germany).}
\address{Krzysztof Bogdan \\ Faculty of Pure and Applied Mathematics,
Wroc\l{}aw University of Science and Technology}
\email{\tt krzysztof.bogdan@pwr.edu.pl}

\author{Tomasz Komorowski}
\thanks{T.K. acknowledges the support of the NCN grant 2020/37/B/ST1/00426.}
 \address{Tomasz Komorowski \\ Institute of Mathematics,
   Polish Academy
  of Sciences\\Warsaw, Poland.}
\email{{\tt tkomorowski@impan.pl}}

 \author{Lorenzo Marino}
 \address{Lorenzo Marino \\ Institute of Mathematics,
   Polish Academy
  of Sciences\\Warsaw, Poland.}
  \email{{\tt lmarino@impan.pl}}


\maketitle

\begin{abstract}
We consider the limit of solutions of scaled linear kinetic
equations with a reflection-transmission-absorption condition at the
interface. Both the coefficient describing the probability of absorption and the scattering
kernel degenerate. We prove that the long-time, large-space limit is the unique solution of a version of the fractional in space heat equation that corresponds to the Kolmogorov equation for a symmetric stable process, which is reflected, or transmitted while crossing the interface and is killed upon the first hitting of the interface.
The results of the paper are related to the work in \cite{kor19}, where the case of a non-degenerate probability of absorption has been considered.

\end{abstract}

\section{Introduction}
We study the asymptotic behaviour of a linear kinetic equation on
the real line $\R$:
\begin{equation}
  \label{eq:8}
  \begin{cases}
 \partial_tW(t,y,k)+ \bar{\omega}'(k) \partial_y W(t,y,k) \, = \,  \gamma L_k W(t,y,k), &\mbox{ on } (0,+\infty)\times\R_\ast\times \T;
 \\
W(0,y,k)\, = \, W_0(y,k), &\mbox{ on } \R\times \T,
  \end{cases}
\end{equation}
with the following boundary condition at $x=0$:
\begin{equation}\label{eq:interface_condition}
\begin{cases}
W(t,0^+, k)\, = \, p_+(k)W(t,0^-,k)+p_-(k)W(t,0^+, -k)+p_0(k)T_o, &\mbox{ on } \T_+;\\
W(t,0^-, k)\, = \, p_+(k)W(t,0^+, k)+p_-(k)W(t,0^-,-k) + p_0(k)T_o, &\mbox{ on } \T_-.
\end{cases}
\end{equation}
Here, $\T$ denotes the unit one-dimensional \textit{torus}, understood as the
interval $[-1/2,1/2]$ with identified endpoints,
$\T_\pm:=[k\in\T:\pm k>0]$ and
$\R_\ast:=\R\setminus\{0\}$. We call the origin $o:=[y=0]$ \textit{interface}.  The parameters $\gamma>0$, $T_o\ge 0$
are given
  and the functions
$\bar{\omega}, p_\pm, p_0$, defined on $\T$, are assumed to be
continuous and non-negative. The
 \textit{scattering operator}
$L_k$, acting only on the variable $k$ in $\T$, is given by
\begin{equation}
\label{L}
L_ku(k)\, := \, \int_{\T}R(k,k')
\left[u(k') - u(k)\right]\, dk'.
\end{equation}
Here $R:\T^2\to[0,+\infty)$ is  $C^2$ smooth and $u$ belongs to $L^1(\T)$.

This equation
arises in the kinetic limit of the evolution of the energy density
in a
stochastically perturbed
harmonic chain interacting with a point Langevin thermostat located at $y=0$,
see \cite{kor19,kors,ko20,ko20a}.
The energy density $W(t,y,k)$ at time $t$ is resolved in both the
\textit{spatial} variable $y\in \R$ and the \textit{frequency} variable $k\in\T$. The function
$\bar\omega$ is the \textit{dispersion relation} of the harmonic chain and
$\bar\om'(k)$ is the \textit{group velocity} of
\textit{phonons} of mode $k$ -- theoretical particles that carry the energy due
to the chain vibrations of frequency $k$.
The presence of
the Langevin thermostat results in the boundary (interface) condition  \eqref{eq:interface_condition}.
The number $T_o$
corresponds to the temperature of the thermostat, while
$p_+(k)$, $p_-(k)$ and  $p_0(k)$ are the respective probabilities of
transmission, reflection and killing of a mode $k$ phonon at
the interface; see \cite{kors}. They are continuous even functions   that satisfy
\begin{equation}
\label{eq:012304}
p_+(k)+p_-(k)+  p_0(k)\, = \, 1, \quad k\in \T.
\end{equation}
With a small risk of ambiguity, in what follows we also write
\begin{equation}
\label{d.ppm0}
p_+:=p_+(0),\quad p_-:=p_-(0),\quad p_0:=p_0(0).
\end{equation}
Our main goal is to study the asymptotic behaviour  of the solutions
of \eqref{eq:8} on macroscopic space-time scales. More precisely, we are
interested in the limits of solutions, when $\lambda\to+\infty$, for the family of \textit{rescaled} equations
\begin{equation}
\label{kinetic-sc0}
\begin{cases}
\dfrac{1}{\lambda}\partial_t
  W_\lambda(t,y,k)+\dfrac{1}{\lambda^{1/\alpha}}\bar\omega'(k)\partial_y W_\lambda(t,y,k)\, = \, \gamma
  L_kW_\lambda(t,y,k), & (t,y,k)\in(0,+\infty)\times\R_\ast\times\T;
\\
W_\lambda(0,y,k)\, = \, W_0(y,k), &(y,k)\in\R\times\T,
\end{cases}
\end{equation}
subject to the boundary condition \eqref{eq:interface_condition},
where
$1<\alpha<2$ is a suitably chosen exponent.

It turns out, see \cite{kob,kor19}, that under appropriate  assumptions on the total scattering kernel $R(k):=\int_{\T}R(k,k')dk'$, group velocity $\bar\om'(k)$ and a suitable choice of $\alpha$,  the limit $\bar W(t,y):=\lim_{\lambda\to+\infty}W_\lambda(t,y,k)$ exists in a weak
(distribution) sense. In the case where \linebreak
$\sup_{k\in\T}[\bar\om'(k)]^2/R(k)<+\infty$ and $\alpha=2$ (diffusive
space-time
scaling), $\bar
W(t,y)$ satisfies the heat equation with the Dirichlet boundary
condition $\bar W(t,0)=T_o$ and $\bar W(0,y)=\int_{\T}W(0,y,k)dk$, see \cite[Theorem 2.2]{kob}.  On the
other hand, when $|\bar\om'(k)|/R(k)\sim |k|^{-\beta}$, for $|k|\ll1$,
with $\beta>1$ and $p_0=\lim_{k\to0}p_0(k)>0$, then letting
$\alpha=1+1/\beta$  one can show, see \cite[Theorem 1.1]{kor19}, that
$\bar W(t,y)$ satisfies a fractional diffusion equation with a
boundary condition at $y=0$. In informal terms, the equation states
that
\begin{equation}
  \label{gen-L}
  \partial_t\bar W(t,y)={\cal L}\Big(\bar W(t,y)-T_o\Big),
\end{equation}
where
${\cal L}$ is the generator of a modified symmetric
$(1+1/\beta)$-stable process. The process behaves like a ``regular'' symmetric
$(1+1/\beta)$-stable process outside the interface $y=0$, but
transmits, reflects, or dies at the interface with the
respective probabilities $p_+$, $p_-$ and $p_0$ of \eqref{d.ppm0}.

In the present paper, we address the situations where $p_0=0$, in fact, we assume that the \textit{absorption probability} $p_0(k)$ satisfies, for some
$\kappa>0$,
\begin{equation}
\label{g}
p_*:=\liminf_{k\to0}\left(\log |k|^{-1}\right)^{\kappa} p_0(k)\in (0,+\infty).
\end{equation}
Furthermore, we suppose that the  \textit{transmission  probability}  does not vanish, that is,
\begin{equation}
\label{p-plus1}
\inf_{k\in\T} p_+(k) >  0 .
\end{equation}
The dispersion relation $\bar\omega:\T\to[0,+\infty)$ is assumed to be even and
\textit{unimodal},  i.e., it possesses exactly one local maximum, at $1/2$, and
one local minimum, at $0$. In addition it is of
the $C^2(\T_*)$ class of regularity, where $\T_*:=\T\setminus\{0\}$, and there exist one-sided limits of its
both derivatives  at
$k=0$. A typical dispersion relation we have in mind  is
$\bar\omega(k)=|\sin(\pi k)|$, which corresponds to the
harmonic chain with the nearest-neighbour interaction.

Concerning the scattering kernel, we assume that it is of the
multiplicative form
\begin{equation}
  \label{R1R2}
R(k,k')\, = \, R_1(k)R_2(k'),
\end{equation} where $R_1$, $R_2$  are two non-negative, even functions  belonging
to  $C^2(\T)$.
Without loss of generality, we also suppose that
\[
\int_{\T}R_j(k)\, dk\, = \, 1, \quad j=1,2.
\]
The functions $R_1$ and  $R_2$ may possibly vanish at some points in
$\T$, but we suppose that there exist exponents $\beta_j>0$, $j=1,2,3$,
such that
\begin{align}
\label{tot}
R^\ast_{j} \, &:= \,\lim_{k\to0}\frac{R_j(k)}{|k|^{\beta_j}}>0,\quad j=1,2,\\
\label{S0}
S_\ast\,&:=\,\lim_{k\to 0 }\frac{ |k|^{\beta_3}|\bar{\omega}'(k)|}{ \gamma R_1(k)}>0,\\
\label{e:Markov}
\beta_1\, &< \, 1+\beta_2,\\
\label{e:alpha}
\alpha\, &:= \, \frac{1+\beta_2}{\beta_3}\, \in \, (1,2).
\end{align}
Consistent with the condition \eqref{e:Markov}, we further require that
\begin{equation}
  \label{Gamma}
{\cal R}:=\int_{\T}\frac{R_2(k)}{ R_1(k)}dk\in (0,+\infty).
\end{equation}
We define the probability measure $\pi(dk):=R_2(k)/\big({\cal R} R_1(k)\big)dk$ and we can \textit{informally} state the main result of the present work; see
 Theorem \ref{main-thm} for a rigorous formulation.
\vspace{10pt} \newline
\textbf{Main Theorem.} \emph{Let $W_0(y,k)$ be a sufficiently regular function satisfying the interface conditions in \eqref{eq:interface_condition} and such that the averaged function}
\begin{equation}
\label{022603-19}
\bar W_0(y)\, := \, \int_{\T}W_0(y,k)\,\pi(dk)
\end{equation}
\emph{satisfies  $\bar{W}_0(0)=T_o$. Let $W_\lambda(t,y,k)$ be the  solution to
\eqref{kinetic-sc0} with interface conditions
\eqref{eq:interface_condition}. Then, under the hypotheses
\eqref{g}--\eqref{e:alpha} made above, the limit
$
\lim_{\lambda\to+\infty}W_\lambda(t,y,k) \, = \, \bar{W}(t,y)$ exists
in a distributional sense on $\R\times \T$ for any $t>0$. Moreover,   $\bar{W}(t,y)$
is a unique weak solution, in the sense of Definition \ref{df011803-19}
below, of the equation:}
\begin{align}
\label{apr204}
&\partial_t \bar{W}(t,y) \,= \, \bar{\gamma}\,  \,{\rm
  p.v.}\int_{yy'>0}q_\alpha (y'-y)[\bar{W}(t,y')-\bar{W}(t,y)] \, dy'\\
&  +\bar{\gamma}p_+\,  \int_{yy'<0}q_\alpha (y'-y)
    [\bar{W}(t,y')-\bar{W}(t,y)]  \, dy'  +\bar{\gamma}p_-\, \int_{yy'<0}q_\alpha (y'-y)[\bar{W}(t,-y')-\bar{W}(t,y)] \, dy',
\nonumber
\end{align}
\emph{with  the initial and  boundary conditions
  $\bar{W}(0,y) = \bar{W}_0(y)$ and
$\bar{W}(t,0)\, = \, T_o$, respectively.
Here, $\bar{\gamma}$ is a positive constant depending on the
parameters of the model (see \eqref{eq:def_bar_gamma} below),
\begin{equation}
\label{cal}
 q_\alpha (y)\, := \, \frac{c_\alpha}{|y|^{ 1+\alpha}}, \qquad
 c_\alpha \, :=\,
 \frac{2^\alpha\Gamma((1+\alpha)/2)}{\sqrt{\pi}\Gamma(\alpha/2-1)},\quad
 y\not=0,
\end{equation}
and $\Gamma(\cdot)$ is the Euler gamma function.}

\begin{remark}
\label{rmk012906-22}
 {\em The above result should be compared with the result of \cite{kor19},
where the case $p_0=\lim_{k\to0}p_0(k)>0$
has been considered. As we have already mentioned, the informal
formulation of the limiting fractional dynamics involves the generator
${\cal L}$ of the process that is symmetric stable and can
be transmitted, reflected or killed when crossing the interface at
$y=0$. In the situation considered in the present paper we can also
view the evolution of $\bar{W}(t,y)$ as described by  equation
\eqref{gen-L}, with ${\cal L}$ corresponding to a symmetric
$\alpha$-stable process that is transmitted or reflected when
crossing the interface at $y=0$, with the respective probabilities
$p_+$ and $p_-$. Furthermore, the process is killed upon the first
hitting of the interface.}
\end{remark}

The fractional diffusion limit of a kinetic equation
has been the subject of intense investigation in recent years. We refer the interested reader to a review of the existing literature contained in
\cite{mellet-rev}. However, there seem to be only few results dealing with a
fractional diffusion limit for kinetic equations with a boundary
condition. In this context, we mention the papers \cite{AS,bgm,ce,ce1,cmp,cmp1,kor19}.
The case that is somewhat related to ours is
considered in \cite{cmp} and \cite{cmp1}. In the first paper, the
convergence of scaled solutions to kinetic
equations in spatial dimension one, with diffusive reflection
condition on the boundary, is investigated.
However, this condition is different from ours. Furthermore, the results of
\cite{cmp}
do not establish
the uniqueness of the limit,
stating only that it satisfies a certain fractional diffusive equation
with a boundary condition and leaving the question of the uniqueness of solutions for the limiting equation open, see the
remark after Theorem 1.2 in~\cite{cmp}.
In the paper \cite{cmp1}, the authors complete the results of
\cite{cmp} and
establish an anomalous
diffusive limit for a family of solutions of  scaled
linear kinetic equations in a one-dimensional bounded domain with
diffusive boundary conditions. The scattering kernel, defined on $\R^2$, is also of
multiplicative form, as in \eqref{R1R2}, and decays according to an appropriate
power law.

Concerning our methods of proof, we rely on the probabilistic interpretation of solutions to kinetic equations. As shown in \cite[Proposition 3.2]{kor19}, the solution can be expressed, with the help of an underlying two-dimensional stochastic process $\{Y(t),K(t)\}_{t\ge 0}$, see Proposition~\ref{prop010701-19a} below. The component $K$ is the frequency mode of a phonon. If outside the interface $[y=0]$, it is described by a
pure jump, $\T$-valued Markov process corresponding to the operator $L_k$ in
\eqref{L}. The component $Y$ of the process is the position of
the phonon. If the phonon does not cross the
interface, it performs a uniform motion between the consecutive scattering events, with velocity $-\bar\om'(K(t))$.
  On the other hand, if
the phonon tries to cross the interface at time $t$, it will be transmitted,
reflected, or killed with the respective probabilities $p_+(K(t))$,
$p_-(K(t))$ and $p_0(K(t))$. In the case of  reflection, the frequency
mode of the phonon changes to $-K(t)$. Using this probabilistic
interpretation, the asymptotic behaviour of the solutions
of \eqref{kinetic-sc0} can be reformulated into the problem of finding
the limit of appropriately scaled processes
$\{Y_\lambda(t),K_\lambda(t)\}_{t\ge0}$, with the parameter $\lambda>0$
corresponding to the ratio between the macro- and microscopic time units. This
approach has been successfully implemented in the case of solutions of
scaled kinetic equations without interface in \cite{bb,jk,jko}. In the
case when the killing probability $p_0$ is strictly positive, as assumed in
\cite{kor19}, one can expect that typically the phonon crosses the interface
only finitely many times before being killed. Its trajectory can then be obtained by a path transformation (consisting in performing suitable reflections)
of the trajectory when no interface is present. This transformation
turns out to be continuous in the Skorokhod $J_1$ topology on the path
space, thus the convergence of the scaled processes in the presence of
the interface is a consequence of the result for the process in the free
space, done, e.g.,  in \cite{jko}. This approach cannot be applied in our present situation, since
with small $p_0$ it is not so obvious how to effectively control the number of interface crossings performed by
the phonon before it is killed. For this reason, we use a different method to deal with the problem. We define a L\'evy-type process
$\big\{Z(t)\big\}_{t\ge 0}$ whose scaled limit $\big\{Z_{\lambda}(t)\big\}_{t\ge0}$, as
$\lambda\to+\infty$, is related to the respective limit of the position of the phonon by a deterministic time change. To find the limit of $\big\{Z_{\lambda}(t)\big\}_{t\ge0}$, we prove that the  associated Dirichlet forms converge, in
the sense of the $\Gamma$-convergence of forms (see Definition
\ref{df2.5} below, or \cite{mosco}), to the
Dirichlet form corresponding to the modified $\alpha$-stable process,
described in Remark \ref{rmk012906-22}. Apart from the limit theory of kinetic equations, we present here some new technical results, which may be of independent interest, for the construction of skew-stable processes, for the characterisation of the fractional Sobolev spaces and for the definition and uniqueness of solutions to non-local parabolic equations.

The paper is organised as follows: after presenting some preliminaries in Section \ref{sec:Basic_notations}, we reformulate the problem of finding the limit of solutions of \eqref{kinetic-sc0}, with the
boundary condition \eqref{eq:interface_condition} into the problem of
finding the limit of the respective stochastic processes. This is done
in Section \ref{sec3}. In Section \ref{sec6.2}, we show how to use the
result concerning the limit of stochastic processes to conclude our
main result, formulated rigorously inTheorem \ref{main-thm}. Section \ref{sec:Proof_converg} is devoted
to the proof of the convergence of the scaled L\'evy
type processes, as stated in Theorem \ref{cor012402-20}, assuming the
$\Gamma$-convergence of their respective Dirichlet forms. The latter is established
in Section \ref{sec4}. Some auxiliary facts are proved in
Appendix \ref{secA} and \ref{appC}.

\setcounter{equation}{0}
\section{Preliminaries}
\label{sec:Basic_notations}

\subsection{Precise assumptions on the model.} Given an arbitrary set $A$ and two functions $f,g:A\to [0,+\infty]$,
we write $f\preceq g$ on $A$ if there is a constant $C>0$ such that
\[
f(a)\, \le \, Cg(a),\quad a\in A.
\]
We also write $f\approx g$ on $A$ when $f\preceq g$ and $g\preceq f$. In what follows, we  always assume the following conventions: $\sum_A f(a) =0$ and $\prod_A f(a)=1$, if $A=\emptyset$, and  $\bar\omega'(0):=0$, even though as a rule $\omega$ is not differentiable at $k=0$.

Let us now denote for any $k$ in $\T_\ast$,
\begin{equation}
\label{bar-t}
\bar{t}(k)\, := \, (\gamma R_1(k))^{-1},
\end{equation}
which can be interpreted as \textit{expected waiting time} for the scattering of a phonon at frequency $k$. The interplay between scattering and \textit{drift} will be captured by the function
\begin{equation}
\label{SS}
S(k)\, := \, \bar{\omega}'(k)\bar{t}(k), \quad k \in [-1/2,1/2].
\end{equation}
Clearly, we can interpret $S(k)$ as the expected distance travelled by the phonon before scattering.
According to our assumptions, $S$ is an odd function,
$C^1 $ smooth on $[-1/2,1/2]\setminus\{0\}$.
To simplify some
considerations below, we also assume that $S$  is a bijection and hence decreases in $(0,1/2]$.  By \eqref{S0},
$\lim_{k\to 0 } |k|^{\beta_3} |S(k)|=S_\ast$ and
we impose some further restrictions on the model coefficients:
\begin{equation}
\label{S01}
\begin{split}
&R_j(k)\, \approx \, |k|^{\beta_j}, \\
&
|S(k)|\, \approx \, \frac{\cos(\pi k)}{|k|^{\beta_3}},\quad |S'(k)|\,
\approx \, \frac{1}{|k|^{1+\beta_3}}\,\, \text{ on }\T_* \quad
\text{and} \quad\\
&
S'_\ast\,:=\,\lim_{k\to0}|k|^{1+\beta_3}|S'(k)| >0.
\end{split}
\end{equation}
 From this point on, we shall assume that all the above hypotheses are in force.

\subsection{Solution of the kinetic equation.} Let us write
$\R_+:=(0, +\infty)$, $\bar{\R}_+:=[0,+\infty)$,  $\R_-:=(-\infty,0)$ and $\bar{\R}_-:=(-\infty, 0]$. Given $T \in \R$, we denote by $\mathcal{C}_{T}$ the class of functions $\phi$ in $C_b(\R_*\times\T_*)$ that can be continuously extended to $\bar{\R}_\iota\times\T_*$, for $\iota\in\{+,-\}$, and satisfy the following interface conditions:
\begin{equation}\label{feb1408_new}
\begin{cases}
\phi(0^+, k)\, = \, p_+(k)\phi(0^-,k)+p_-(k)\phi(0^+, -k)+p_0(k)T, &\mbox{ on } \T_+;\\
\phi(0^-, k)\, = \, p_+(k)\phi(0^+, k)+p_-(k)\phi(0^-,-k) + p_0(k)T, &\mbox{ on } \T_-.
\end{cases}
\end{equation}
Clearly, the constant function $T_o$ belongs to $\mathcal{C}_{T_o}$. Furthermore,
$F\in\mathcal{ C}_{T_o}$ if and only if $F-T_o\in\mathcal{C}_0$.
This allows us to reduce the proofs of some results below to the case $T_o=0$. Following \cite{kor19}, we now recall the definition of solution to equation \eqref{eq:8} with the interface conditions \eqref{eq:interface_condition}.
\begin{df}
\label{df013001-19}
A bounded, continuous function $W\colon
\R_+\times\R_*\times \T_*\to \R$ is called a solution to equation \eqref{eq:8} with the interface conditions \eqref{eq:interface_condition} if all the following conditions are satisfied:
\begin{itemize}
    \item for any $\iota, \iota'$ in $\{-,+\}$, the restriction of $W$ to $\R_+\times\R_\iota\times \T_{\iota'}$ can be extended to a bounded, continuous function on
$\bar\R_+\times\bar\R_{\iota}\times \bar\T_{\iota'}$;
\item for any $(t,y,k)\in \R_+\times\R_*\times  \T_*$ fixed, the function $s\mapsto W(t+s,y+\bar\omega'(k) s,k)$ is continuously differentiable in a neighbourhood of $s=0$, the directional derivative
\begin{equation}
\label{Dt}
D_tW(t,y,k)
:=\frac{d}{ds}_{|s=0} W(t+s,y+\bar\omega'(k) s,k)
\end{equation}
is bounded in $\R_+\times\R_*\times \T_*$, and
\[
 D_tW(t,y,k) \, = \, \gamma L_k W(t,y,k),
\quad (t,y,k)\in\R_+\times \R_*\times\T_*;
\]
\item the interface conditions in \eqref{eq:interface_condition} hold, together with the initial condition
\[
\lim_{t\to0+} W(t,y,k)\, = \, W_0(y,k),\quad (y,k)\in\R_*\times\T_*.
\]
\end{itemize}
\end{df}
We highlight that, at least formally, the directional derivative satisfies
\[ D_tW(t,y,k)\, = \, \left[\partial_t+\bar\omega'(k)\partial_y\right]W(t,y,k),\]
which justifies the notation in \eqref{Dt}. It has been shown in \cite{kor19} that there exists a unique (classical) solution, in the sense of Definition \ref{df013001-19}, to the Cauchy problem \eqref{eq:8}, if  the initial distribution $W_0$ belongs to $\mathcal{C}_{T_o}$. On this regard, see as well Proposition \ref{prop010701-19a} below.

\subsection{Weak solution of the limit equation}
  \label{ws}

Recalling the definitions of $p_\pm$ in \eqref{p-plus1}, we  introduce the \textit{bilinear form}
\begin{multline}
\label{hat-C1}
\hat{\cal E}[u,v]\, := \, \frac{1}{2} \int_{\R^2}(u(y')-u(y))(v(y')-v(y)) q_\alpha (y'-y)\left(\mathds{1}_{yy'>0}
    +p_+\mathds{1}_{yy'<0}\right)dydy'\\
+
\frac12\int_{\R^2}(u(-y')-u(y))(v(-y')-v(y)) q_\alpha (y'-y)p_- \mathds{1}_{yy'<0}\,dydy',
\end{multline}
and the associated \textit{quadratic form} $\hat{\cal E}[u]:=\hat{\cal E}[u,u]$.
We can now define the semi-norm $\|u\|_{\mathcal{H}_o}\, := \, \hat{\cal E}^{1/2}[u]$ for any Borel function $u\colon \R_\ast \to \R$ such that the expression is finite and then denote by $\mathcal{H}_o$ the completion of $C_c^\infty(\R_*)$ under $\Vert \cdot \Vert_{\mathcal{H}_o}$.
We proceed to the definition of   a weak solution to   \eqref{apr204}  with the interface conditions \eqref{eq:interface_condition}.
\begin{df}
\label{df011803-19}
A bounded  function
$\bar{W}\colon \bar{\R}_+\times \R\to \R$ is  a  weak
solution to equation \eqref{apr204} with interface conditions \eqref{eq:interface_condition}
if
\begin{itemize}
\item[i)] $\bar{W}(\cdot)-T_o\in L^2_{\text{loc}}([0,+\infty);\mathcal{H}_o)$ and $\bar{W}(\cdot)-T_\infty\in C([0,+\infty);L^2(\R))$ for some $T_\infty\in\R$,
\item[ii)]
  for any  $F\in C^\infty_c([0,+\infty)\times\R_*)$ and any $t\ge 0$, it holds that
\begin{multline}
\label{021803-19xx}
\int_{\R}F(0,y)[ \bar{W}_0(y)-T_o]\, dy
  =\int_{\R}F(t,y)[\bar{W}(t,y)-T_o]\, dy\\
-\int_0^{t}\int_{\R}\partial_sF(s,y) [\bar{W}(s,y)-T_o]\, dy ds +\bar{\gamma} \int_0^{t}\hat{\cal
  E} [F(s,\cdot),\bar{W}(s,\cdot)-T_o] \, ds.
 \end{multline}
 \end{itemize}
\end{df}
By Proposition \ref{prop:prob_repres_limit_sol} below, our definition
guarantees the uniqueness of solutions.

\subsection{Statement of the main result.} Once we have formulated the
definition of a solution, we are ready to formulate rigorously  our main result. It reads as follows:
\begin{thm}
\label{main-thm}
Suppose that the assumptions \eqref{g}, \eqref{p-plus1} and \eqref{S01}
are in force. Let $T_\infty\in \R$ and $W_0$ in ${\cal C}_{T_o}$ be such that $\bar{W}_0-T_\infty \in L^2(\R)$ and $\bar W_0-T_o\in \mathcal{H}_o$, with $\bar{W}_0$ defined in \eqref{022603-19}. Let $W_\lambda(t,y,k)$ be the solution, in the sense of Definition \ref{df013001-19}, to the Cauchy problem
\eqref{kinetic-sc0}. Then,
\begin{equation}
\label{conv}
\lim_{\lambda\to+\infty}\int_{\R\times \T}W_\lambda(t,y,k)F(y,k)\, dkdy
\, = \, \int_{\R\times \T}\bar{W}(t,y)F(y,k)\, dkdy,
\end{equation}
for all
$t>0$ and test functions $F\in C^\infty_c(\R\times\T)$. Moreover, the limit $\bar{W}(t,y)$
is the weak solution, in the sense of Definition $\ref{df011803-19}$, to the Cauchy problem \eqref{apr204}
with the initial distribution $\bar W_0$ and the fractional diffusion coefficient
\begin{equation}
\label{eq:def_bar_gamma}
    \bar{\gamma}\, := \, \gamma R^\ast_2S^{1+\alpha}_0 \Gamma(\alpha+1)/S'_\ast.
\end{equation}
\end{thm}

\setcounter{equation}{0}

\section{Probabilistic representation of solutions}\label{sec3}

\subsection{Construction of the position/momentum processes} As usual, $\N:=\{1,2,\ldots\}$ and $\N_0:=\{0,1,2,\ldots\}$. Let $(\Omega,\mathcal{F},\mathbb{P})$ be a probability space carrying the following random objects. We
consider
a Markov chain $K_n(k)$
reporting the consecutive frequencies of the phonon,
and the renewal process $\fT_n(k)$ of the
scattering times of the phonon frequency (this is Fraktur font $T$).
More precisely, $\{K_n(k)\}_{n\in \N_0}$ is such that $K_0(k)=k$ and $\{K_n(k)\}_{n\in \N}$ are i.i.d.\ random variables on $\T$, distributed
according to the following probability measure:
\begin{equation}\label{barR}
\mu(dk)\, := \, R_2(k)dk,
\end{equation}
and
\begin{equation}
\label{010502-19}
\fT_n(k)\, := \,
\sum_{0\le j<n}\bar t(K_j(k))\tau_j, \quad n \in \N_0,
\end{equation}
where $\{\tau_n\}_{n\in \N_0}$ is an independent sequence of i.i.d.\ exponentially distributed random variables with intensity $1$ and, we recall, $\bar{t}(k)$ was defined in \eqref{bar-t}.
We introduce the process $\{\fT(t,k)\}_{t\ge0}$ as the linear interpolation between the values of $\fT_n(k)$:
\[
\fT(t,k)=\fT_n(k)+(t-n)(\fT_{n+1}(k)-\fT_n(k))\quad \mbox{if}\quad  t\in [n, n+1),\;n\in \N_0.
\]
We then define the continuous-time frequency (momentum) process as
\[K(t,k)\, := \, K_{[\fT^{-1}(t,k)]}(k)=K_n(k) \quad \mbox{if} \quad t\in [\fT_n(k), \fT_{n+1}(k)),\; n\in \N_0.\]
Here $\fT^{-1}$ is the inverse function of $t\mapsto \fT(t,k)$ and $[\cdot]$ denotes the integer part.

According to \eqref{SS}, the phonon position at the time of the $n$-th scattering of its frequency is
\[Z_n(y,k)\, :=\,  y- \sum_{0\le j<n} S(K_j)\tau_j, \quad n \in \N_0.\]
In particular,
the law of $Z_n(y,k)$, for each $n\in \N$, is absolutely continuous with respect to the
Lebesgue measure on $\R$.
We then consider an auxiliary Poisson process
$\{N(t)\}_{t\ge0}$ of intensity $1$ that is independent of both $\{K_n\}_{n\ge0}$ and $\{\tau_n\}_{n\ge0}$. We define
$\{\tilde N(t)\}_{t\ge0}$ as the linear
interpolation between the nodal points of $N$. Namely, if $n:=N(t)$ then
\[
\tilde N(t):=n+(t-l)/(r-l),
\]
where $l:=\inf\{s:N(s)=n\}$ and $r:=\inf\{s: N(s)=n+1\}$.
Our next process, $\{\tilde Z(t,y,k)\}_{t\ge0}$, is obtained by linearly interpolating between the nodal points of $Z_{N(t)}(y,k)$: if $n=N(t)$, then
\[
\tilde Z(t,y,k):=Z_n(y,k)+(Z_{n+1}(y,k)-Z_n(y,k))(t-l)/(r-l).
\]
We next define the transmission/reflection/absorption mechanism at
the interface $o=[y=0]$. To this end, we fix $(y,k)\in \R_*\times
\T_*$ and consider the times
$\left\{n_{m}\left(\Big\{Z_n(y,k)\Big\}_{n\ge0}\right)\right\}_{m\ge0}$
when $\Big\{Z_n(y,k)\Big\}_{n\ge0}$  crosses the interface.   Namely, we let $n_0:=0$ and then define recursively
\begin{equation}
\label{eq:def_stopping_times_n}
n_{m+1} \, := \, \inf\left\{n>n_{m}\colon   (-1)^m y Z_n(y,k)<0\right\}, \quad m =0,1,\ldots.
\end{equation}
For the sake of brevity, when there is no danger of confusion, we skip the sequence $\Big\{Z_n(y,k)\Big\}_{n\ge0}$
from the notation.
 Similarly, we define the sequence $\Big\{\tilde{\mathfrak{s}}_{m}\Big(\{\tilde Z(t,y,k)\}_{t\ge0} \Big)\Big\}_{m\ge0}$ of consecutive times when the process  $\tilde Z(t,y,k)$ crosses the interface $o$. Namely, we let $\tilde{\mathfrak{s}}_{0}=0$ and
\begin{equation}
\label{eq:def_stopping_times_frak_s}
\tilde{\mathfrak{s}}_{m+1}\, := \, \inf\left\{t>\tilde{\mathfrak{s}}_{m}\colon   (-1)^m y\tilde{Z}(t,y,k)<0\right\}, \quad m =0,1,\ldots.
\end{equation}
Again, we simplify the notation by omitting the path of the process when
there is no danger of confusion. Sometimes, to highlight the
dependence of the crossing times on the starting point (which will be
relevant in the argument), we may write
$\tilde{\mathfrak{s}}_{m,y,k}$, or $\tilde{\mathfrak{s}}_{m,y}$.

We define ${\cal S}(t,k):=\tilde N^{-1}\left(\fT^{-1}(t,k)\right)$ and ${\mathfrak{s}}_{m}:=\mathcal{S}(\tilde{\mathfrak{s}}_{m},k)$. Then, $n_{m}=N(\mathfrak{s}_{m})$ and
\[
\tilde{\mathfrak{s}}_{m} \, < \, \mathfrak{s}_{m} \, < \, \tilde{\mathfrak{s}}_{m+1},\quad \mathbb{P}\text{-a.s.}
\]
We now consider a sequence $\{\sigma_{m}\}_{m\ge0}$ of
$\{-1,0,1\}$-valued random variables, that are independent when
conditioned on $\{K_n(k)\}_{n\ge0}$, such that $\sigma_0:=1$ and
\begin{equation}\label{d.rim}
\bbP\left(\sigma_{m}=\iota|\{K_n(k)\}_{n\ge0}\right) \,= \, p_\iota(K_{n_{m}-1}(k)), \quad \iota \in \{-1,0,1\},\,m\in \N.
\end{equation}
Here and below, $p_\iota$ means $p_\pm$ if $\iota=\pm1$, respectively.
Of course, $\{\sigma_m\}_{m\ge 1}$ can be \textit{defined} by applying the quantile functions for \eqref{d.rim} to independent i.i.d. uniform random variables.
We can finally add the random interface mechanism to the processes considered. Namely,
\begin{equation}
\label{eq:def_interface_process}
\tilde{Z}^o(t,y,k) \,:= \,\left(\prod_{j=1}^{m}\sigma_{j}\right) \tilde{Z}(t,y,k), \quad
K^o(t,k) \,:= \,\left(\prod_{j=1}^{m}\sigma_{j}\right) K(t,k),
\quad t \in [\tilde{\mathfrak{s}}_{m},\tilde{\mathfrak{s}}_{m+1}),
\end{equation}
for $m\in \N_0$.
The ``true'' interface position process is then given by
\begin{equation}
\label{011404-20}
Y^o(t,y,k)\, := \, \tilde
Z^o({\cal S}(t,k),y,k) \, = \, y-\int_0^{t}\bar\omega'(K^o(s,k))\, ds.
\end{equation}
We now denote by ${\mathfrak f} := \min\left\{m\in \N\colon \sigma_{m}=0\right\}$ the interface crossing at which the particle gets absorbed and let $\mathfrak{ s}_{{\mathfrak f}} :={\cal S}(\tilde{\mathfrak
s}_{{\mathfrak f}},k)$. The following probabilistic representation of
the solution to   \eqref{eq:8}, with interface conditions
\eqref{eq:interface_condition}, holds.
\begin{prop}
\label{prop010701-19a}
Let $W_0$ be in $\mathcal{C}_{T_o}$ such that $W_0(0,k)=T_o$ on $\T_\ast$. Then, the function
\begin{equation}
\label{010702-20a}
 W(t,y,k)\, = \, \mathbb{E}\left[W_0\left(Y^o(t,y,k),K^o(t,k)\right)\right], \qquad (t,y,k)\in\bar \R_+[0,+\infty)\times \R_\ast\times \T_\ast,
\end{equation}
is the unique classical solution, in the sense of Definition \ref{df013001-19}, to the Cauchy problem \eqref{eq:8} with initial distribution $W_0$.
\end{prop}
\proof
The existence and uniqueness of a classic solution and the following representation
\begin{equation}
\label{010701-19a}
W(t,y,k)\, = \, \mathbb{E}\left[W_0\left(Y^o( t,y,k), K^o( t,k)\right),\,
  t< \mathfrak{ s}_{y,k,{\mathfrak f}}\right]+T_o\mathbb{P}\left[
  t\ge \mathfrak{ s}_{y,k,{\mathfrak f}}\right],
\end{equation}
have already been shown in
\cite{kor19}, Section $A$, and Proposition $3.2$. Recalling that $Y^o( t,y,k)=0$ if $t\ge {\mathfrak{s}}_{y,k,{\mathfrak f}}$, we can use that $W_0(0,k)=T_o$ to get that
\[T_o\mathbb{P}\left[
  t\ge \mathfrak s_{y,k,{\mathfrak f}}\right] \, = \,\mathbb{E}\left[W_0\left(0, K^o( t,k)\right),\,
  t\ge \mathfrak{ s}_{y,k,{\mathfrak f}}\right]\, = \, \mathbb{E}\left[W_0\left(Y^o( t,y,k), K^o( t,k)\right),\,
  t\ge \mathfrak{ s}_{y,k,{\mathfrak f}}\right].\]
Equation \eqref{010702-20a} then follows immediately from \eqref{010701-19a}.
\qed

\subsection{Scaling of the process $(Y^o( t,y,k), K^o( t,k))$}
\label{sec4.3}
Equations \eqref{kinetic-sc0} may be considered a special case of \eqref{eq:8}, so we first focus on establishing flexible notation for later use.
As before, we fix $(y,k)\in\R_*\times \T_*$. We let $\lambda>0$ and \textit{rescale} the \textit{clock processes}, introduced in the previous section, as follows:
\[
N_\lambda(t) \, := \, N(\lambda t), \quad \tilde{N}_\lambda(t) \, := \, \tilde{N}(\lambda t), \quad
\fT_\lambda(t,k) \, := \, \frac{1}{\lambda} \fT(t,k)
\]
and
\begin{equation}
\label{010502-19N1}
{\cal S}_\lambda(t,k)\, := \, \frac{1}{\lambda}{\cal S}(\lambda t,k)
\, = \,
\tilde N_{\lambda}^{-1}\left(\fT_\lambda^{-1}(t,k)\right).
\end{equation}
We then notice that for \textit{large} $\lambda$, the process ${\cal
  S}_\lambda(t,k)$ becomes almost deterministic. More precisely, by \eqref{Gamma} we have that
\begin{equation}
\label{bar-theta}
\bar\theta\, := \, \left(\bbE [\bar t(K_1)]\right)^{-1}\, = \,
\gamma{\cal R}^{-1} \, < \, +\infty.
\end{equation}
Except for the first
deterministic term, the inverse ${\cal S}_\lambda^{-1}(s,k) =
\fT_\lambda(\tilde N_{\lambda}(s),k)$ can be represented as a Poisson
sum of i.i.d.\ variables with finite first moment and so, by a standard argument using the strong  law of large numbers, it can be shown that:
\begin{prop}
\label{prop011502-20}
For any $t_*\in [0,+\infty)$, we have
\begin{equation}
\label{021502-20}
\lim_{\lambda\to+\infty}\sup_{t\in[0,t_*]}\left|{\cal S}_\lambda(t,k)-\bar{\theta}t\right|=0,\quad \bbP\mbox{-a.s.}
\end{equation}
\end{prop}
We can now rescale the ``position'' process in the following way. Let us define $\{\tilde Z_\lambda(t,y,k)\}_{t\ge0}$ as the linear interpolation between the nodal points of  $Z_\lambda(t,y,k)$, where
\begin{align}
\label{eq:def_Z_lambda}
Z_\lambda(t,y,k)\, &:= \, Z_{N_\lambda(t)}^{\lambda}(y,k),\\
\label{eq:def_Z_n_lambda}
Z_{n}^{\lambda}(y,k)\, &:= \, y-   \frac{1}{\lambda^{1/\alpha }}\sum_{j=0}^{n-1}
  S(K_j)\tau_j, \quad n \in \N_0.
\end{align}
We then construct the sequences $\{{\bar{\mathfrak{s}}}_{m}^{\lambda}\}_{m\ge0}$,
$\{{\tilde{\mathfrak{s}}}_{m}^{\lambda}\}_{m\ge0}$ of stopping times for the rescaled processes $Z_\lambda(t,y,k)$, $\tilde{Z}_\lambda(t,y,k)$ as in \eqref{eq:def_stopping_times_frak_s}, skipping $(y,k)$ from the notation, when unambiguous. Since the law of $ Z_\lambda(t,y,k)$ is absolutely continuous, for each $y\in\R$ there exists a strictly increasing sequence
$n_{m}^\lambda:=N_\lambda({\bar{\mathfrak{s}}}_{m}^{\lambda})$ such that
\[
{\bar{\mathfrak{s}}}_{m}^{\lambda} \, = \, \tilde N_\lambda^{-1}\left(n_{m}^\lambda\right).
\]
To set the notation for the processes subject to the random mechanism at the interface, let $\{\sigma_{m}^\lambda\}_{m\in \N}$ be a sequence of $\{-1,0,1\}$-valued random variables that are independent when conditioned on $\{K_n(k)\}_{n\ge0}$ and such that
\[
\bbP\left(\sigma_{m}^\lambda=\iota|\{K_n(k)\}_{n\ge0}\right) \, = \, p_\iota(K_{n_{m}^\lambda-1}(k)),\qquad
\iota\in \{-1,0,1\}.\]
We then set
\[
{\mathfrak f}^\lambda\, := \, \min\{m \ge1\colon \sigma_{m}^{\lambda}=0\}, \qquad \tilde{\mathfrak{s}}_{\mathfrak f}^\lambda\, := \, \tilde{\mathfrak{s}}^{\lambda}_{{\mathfrak{f}^\lambda}}, \qquad  {\bar{\mathfrak{s}}}_{\mathfrak {f}}^\lambda \, := \, \bar{\mathfrak{s}}^{\lambda}_{{\mathfrak f}^\lambda}.\]
The jump process $\{Z_\lambda^o(t,y,k)\}_{t\ge0}$ is now defined, for any $m\in \N_0$, as
\begin{equation}
    \label{eq:def_Z_o_lambda}
Z_\lambda^o(t,y,k)\, := \, \left(\prod_{j=0}^{m}\sigma_{j}^{\lambda}\right) Z_\lambda(t,y,k),\quad
 t\in[\bar{\mathfrak{s}}^\lambda_{m}, \bar{\mathfrak{s}}^\lambda_{m+1}).
\end{equation}
Here $\sigma_{0}^{\lambda}:=1$.
Similarly, the continuous trajectory process $\{\tilde
Z_\lambda^o(t,y,k)\}_{t\ge0}$ can be obtained as in \eqref{eq:def_Z_o_lambda} but with respect to the stopping times $\{\tilde{\mathfrak{s}}^{\lambda}_m\}_{m\in\N_0}$.

In what follows, we show that it is unlikely that the scaled position
process $Z_\lambda(t,y,k)$ crosses the interface after the first
jump, when $\lambda$ becomes large. Indeed, recall from
\eqref{eq:def_Z_n_lambda} that
$Z_1^\lambda(y,k)=y-\lambda^{-1/\alpha}S(k)\tau_0$ where $\tau_0$ is
$\exp(1)$-distributed. Let
\begin{equation}
\label{Ala}
A^{\la}(y,k):=\left[ yZ_1^\lambda(y,k)<0\right]
\end{equation}
 The following estimate follows immediately.
\begin{prop}
\label{prop012002-20}
For all $\lambda>0$, $y\in\R_*$ and $k\in \T_*$ we have
$\bbP\left( A^{\la}(y,k)\right)\, \le \,
\exp\left\{-\frac{|y|\lambda^{1/\alpha }}{|S(k)|}\right\}.$
\end{prop}

\subsection{Auxiliary stable L\'evy process}
\label{AsLp}

To describe the limit, as $\lambda$ goes to $+\infty$, of processes $\{Z_\lambda^o(t)\}_{t\ge 0}$, we consider the $\alpha$-stable L\'evy process $\{\zeta(t)\}_{t\ge 0}$ with the infinitesimal generator
\begin{equation}
\label{gen-stable}
\mathcal{L}_{\alpha}u(y)\, := \, \bar{r}_\ast \,{\rm
  p.v.}\int_{\R}[u(y')-u(y)]q_\alpha(y'-y)\, dy',\quad u\in C_c^2(\R),
\end{equation}
where $\bar{r}_\ast \, := \,   R^\ast_2S_\ast^{1+\alpha} \Gamma(\alpha+1)/S'_\ast$.
For any $y\neq 0$, we denote $\zeta(t,y) := y+\zeta(t)$. By a~straightforward modification of  \eqref{eq:def_stopping_times_frak_s},
we let $\mathfrak{t}_{0}:=0$ and denote by $\Big\{\mathfrak{t}_{m}\Big(\Big\{\zeta(t,y)\Big\}_{t\ge0}\Big)\Big\}_{m\ge1}$ the consecutive times when the process $\zeta(t,y)$ crosses the interface
$o$. We often abbreviate
$\mathfrak{t}_{m,y}:=\mathfrak{t}_{m}\Big(\Big\{\zeta(t,y)\Big\}_{t\ge0}\Big)
$ and we also let
\[
\mathfrak{t}_{y,\mathfrak{f}} \,:=\, \inf\left\{t>0\colon \zeta(t,y)=0\right\}.\]
It is known that $\mathfrak{t}_{y,\mathfrak{f}}$ is finite
$\mathbb{P}$-a.s.\ and its law is absolutely continuous with respect
to the Lebesgue measure on $\R$ (cf.\ \cite[Example
$43.22$]{sato}).
We can finally introduce the tentative limit process
$\{\zeta^o(t,y)\}_{t\ge 0}$, as follows.
Recalling the definitions of $p_{\pm}$ in \eqref{d.ppm0} and $p_0=0$, stemming from \eqref{g}, we take a sequence $\{\sigma_m\}_{m\in \N}$ of i.i.d.\ $\{-1,1\}$-valued random
variables such that for any $m\in \N$, the random variable
$\sigma_m$ is independent of  $\{\zeta(t,y)\}_{t\ge0}$ and
$\mathbb{P}(\sigma_m=\pm1) = p_\pm$ (the notation of \eqref{d.rim} is best ignored from now on). Letting $\sigma_0:=1$, we define for $m\in \N_0$,
\begin{equation}
\label{Zkr}
\zeta^o(t,y)\, := \,\left(\prod_{j=0}^{m}\sigma_j\right)
\zeta(t,y), \quad \mbox{if} \quad t \in[\mathfrak{t}_{y,m},\mathfrak{t}_{y,m+1}),
\end{equation}
and $\zeta^o(t,y):= 0$ if $t\ge \mathfrak{t}_{y,\mathfrak{f}}$. The process, discussed in detail below, is an interesting take on the question of construction of \textit{skew} stable L\'evy processes, a topic recently discussed in \cite{2021-AI-AP-a}. However, in view of limit theorems for kinetic equations, the following is our main motivation.
\begin{thm}
\label{cor012402-20}
When $\lambda$ tends to $+\infty$,
the processes $\left\{Z_{\lambda}^o(t,y,k)\right\}_{t\ge0}$
converge both in finite distributions and weakly, over $\mathcal{D}[0,+\infty)$ with the $J_1$-topology, to  $\left\{\zeta^o(t,y)\right\}_{t\ge0}$.
\end{thm}
The proof of the  theorem  shall be presented in Section \ref{sec:Proof_converg} below.

For simplicity, denote by $\mathcal{X}$ the space $\mathcal{D}[0,+\infty)\times \mathcal{C}[0,+\infty)$ equipped with the product of the $M_1$-topology and the topology of uniform convergence over compact intervals.
The definition of the $M_1$-topology on $\mathcal{D}[0,+\infty)$ can be found, e.g., in \cite[Section 12]{whitt}. It then immediately follows from Theorem \ref{cor012402-20}, Proposition \ref{prop011502-20}
that:
\begin{cor}
\label{cor020203-20}
As $\lambda \to +\infty$, the processes $\left\{\left(\tilde Z_{\lambda}^o(t,y,k),\mathcal{S}_{\lambda}(t,k)\right)\right\}_{t\ge0}$ converge both in finite distributions and weakly, over $\mathcal{X}$, to $\left\{\left(\zeta^o(t,y),\bar{\theta}t\right)\right\}_{t\ge0}$.
\end{cor}

Let us denote $K_\lambda(t,k):=K(\lambda t,k)$ and
$
\mathfrak{s}_{m}^\lambda:={\cal S}_{\lambda}(\tilde {\mathfrak{s}}_{m}^\lambda,k)
$. We finally define the position- momentum  process
$\Big(K_\lambda^o(t,k), Y_\lambda^o(t,k,y)\Big)_{t\ge0}$   as
\begin{align}
K_\lambda^o(t,k)\, &:= \, \left(\prod_{j=1}^{m}\sigma_{j}^{\lambda}\right) K_\lambda(t,k),\quad
 t\in[\mathfrak{s}_{m}^\lambda, \mathfrak{s}_{m+1}^\lambda), \quad m \in \N_0 ;
\\
\label{011502-20}
Y^o_\lambda(t,y,k) \, &:= \, \tilde Z^o_\lambda({\cal S}_\lambda(t,k),y,k)\, = \, y-\frac{1}{\lambda^{1-1/\alpha}}\int_0^t\bar\omega'( K^o_\lambda(t,s)) \, ds.
\end{align}
If we   let $\eta^o(t,y) := \zeta^o\left(\bar{\theta}t,y\right)$, where $\bar{\theta}$ is defined in \eqref{bar-theta}, then
the following result is an immediate consequence of Corollary \ref{cor020203-20}.
\begin{cor}
\label{cor030203-20}
As $\lambda$ goes to $+\infty$, the processes $\left\{Y^o_{\lambda}(t,y,k)\right\}_{t\ge0}$
converge both in finite distributions and weakly, over $\mathcal{D}[0,+\infty)$ with the $M_1$-topology, to the process
$\left\{\eta^o(t,y)\right\}_{t\ge0}$.
\end{cor}
\proof
Invoking Theorem $7.2.3$ of \cite{whitt-s} and using formula  \eqref{011502-20}, we conclude
the weak convergence of the processes. Since,
for any deterministic time $t\ge 0$ we have
\[\bbP[\eta^o(t,y)=\eta^o(t-,y)]\, = \, 1,\]
by virtue of Lemma $6.5.1$ in \cite{whitt}, the set of discontinuities of the one-dimensional projection mapping
$\omega\mapsto \eta^o(t,y,\omega)$ is of null probability. Using the
continuous mapping theorem, see Theorem $2.7$ of \cite{billingsley}, we
conclude the convergence of the one-dimensional distributions. The
generalisation to finite-dimensional distributions is trivial.
\qed

We conclude this section by presenting a probabilistic characterisation of the weak solution of the limit Cauchy problem \eqref{apr204} in terms of the process $\eta^o(t,y)$. A proof of this result can be found in the Appendix \ref{appC} below.

\begin{prop}
\label{prop:prob_repres_limit_sol}
Let $\bar{W}_0-T_o$ be in $\mathcal{H}_0$ such that $\bar{W}_0-T_\infty$ is in $L^2(\R)$ for some $T_\infty\in \R$.
Then,  the function
\begin{equation}
\label{eq:prob_repr_bar_W}
\bar{W}(t,y) \, = \, \bbE\left[\bar W_0\left(\eta^o(t,y)\right)\right], \qquad (t,y) \in \bar\R_+\times \R_\ast,
\end{equation}
is the unique weak solution, in the sense of Definition \ref{df011803-19}, to the Cauchy problem \eqref{apr204} with initial condition $\bar{W}_0$.
\end{prop}

\section{Proof of Theorem \ref{main-thm}}
\label{sec6.2}

Since $W_{\lambda}(t,y,k)-T_o$ is the solution  of
\eqref{kinetic-sc0} with the interface condition
\eqref{eq:interface_condition} corresponding to the zero thermostat temperature, we
assume without loss of generality that $T_o=0$.
Clearly, the solution $W_\lambda$ is given by \eqref{010702-20a}. In what follows, we often write $W_\lambda(t)$ for $W_\lambda(t,\cdot,\cdot\cdot)$.

 Fix   a test function $F$ in $C^\infty_c(\R\times
 \T)$. Suppose that ${\rm supp}\,F\subset [-M_0,M_0]\times\T$.
 Thanks to Corollary \ref{cor030203-20} for any $\eps>0$  we can find a
sufficiently large $M>1$ such that
$$
\limsup_{\lambda\to+\infty}\bbP[|Y_\lambda^o(t,y,k)|\ge M]<\eps,\quad
(y,k)\in [-M_0,M_0]\times\T.
$$
This fact allows us to restrict, without loss of generality, our
attention to the case
where $T_\infty=0$.
By a standard approximation argument, it is further enough to consider an
initial condition $W_0\in C^\infty_c(\R\times\T)\cap
\mathcal{C}_{0}$, i.e. in particular
$W_0(0,k)=0$, for any $k$ in $\T$.
If the initial condition $W_0$ is independent of $k$, i.e. $W_0(y,k)=W_0(y)$, then we can immediately conclude the proof from Proposition \ref{prop010701-19a} (with respect to the rescaled processes), equation \eqref{eq:prob_repr_bar_W} and Corollary \ref{cor030203-20}.

The next
result enables  to
replace  an arbitrary initial condition $W_0(y,k)$ by its
average over $\T$ with respect to the measure
$\pi(dk)=R_2(k)/\big({\cal R}R_1(k)\big)dk$. For notational
simplicity, we denote by $L^2_{\pi}(\R\times\T)$ the $L^2$-space on $\R\times\T$ with respect to the product measure $\pi(dk)dy$.
\begin{lemma}
\label{lm010203-20}
Suppose that $W_0\in C_c(\R\times\T) $.
Let $\bar W_0$ be the function defined by \eqref{022603-19}. Then, for any $\eps>0$, there exists ${\delta}_0>0$ such that
\begin{equation}
\label{012201-20}
\limsup_{\lambda\to+\infty}\| W_\lambda(\delta)-\bar W_0\|_{L^2_{\pi}(\R\times\T)}\, < \, \eps, \qquad \forall \, \delta \in [0,{\delta}_0).
\end{equation}
\end{lemma}
The above lemma shall be proved at the end of this
section. Furthermore, computing $\frac{d}{dt}
\|W_{\lambda}(t)\|^2_{L^2_\pi(\R\times\T)}$ and using the equation
\eqref{kinetic-sc0} together with the interface conditions
\eqref{eq:interface_condition}, see the  calculations in \cite[Section
3]{kob}, we immediately conclude the following.
\begin{prop}
\label{cor010203-20}
The norm $  \|W_{\lambda}(t)\|_{L^2_\pi(\R\times\T)}$, defined as a
function of
$t\in [0,\infty)$, does not increase.
\end{prop}
We can now finish the proof of  Theorem \ref{main-thm}. We first show the following variant of \eqref{conv},
\begin{equation}
\label{conv1}
\lim_{\lambda\to+\infty}\int_{\R\times \T}W_\lambda(t,y,k)F(y,k)\, \pi(dk)dy
\, = \, \int_{\R\times \T}\bar{W}(t,y)F(y,k)\, \pi(dk)dy.
\end{equation}
Choose an arbitrary $\eps>0$. Let $\delta >0$ be as in Lemma \ref{lm010203-20} and $\tilde{W}_\lambda$ the solution to the Cauchy problem \eqref{kinetic-sc0} when the initial
condition is given at time $t=\delta$ by $\tilde{W}_\lambda(\delta,y,k)=\bar W_0(y)$. In particular, we can use Proposition \ref{prop010701-19a} to represent $\tilde{W}_\lambda$ as
\begin{equation}
\label{eq:prob_repr_tilde_W_la}
\tilde{W}_\lambda(t,y,k)\, = \, \bbE\left[\bar W_0\left(Y^o_\lambda( t-\delta,y,k)\right)\right], \qquad t \ge \delta.
\end{equation}
By virtue of Proposition
\ref{cor010203-20} and Lemma \ref{lm010203-20}, we have that
\[
\limsup_{\lambda\to +\infty}\| W_\lambda(t)-\tilde{W}_\lambda(t)\|_{L^2_{\pi}(\R\times\T)}\, \le \, \limsup_{\lambda\to +\infty}\| W_\lambda(\delta)-\bar W_0\|_{L^2_{\pi}(\R\times\T)}\, < \, \eps.
\]
Therefore, it follows that
\[
\limsup_{\lambda\to+\infty}\left|\int_{\R\times\T}F(y,k)\left[
  W_\lambda(t,y,k)- \tilde{W}_\lambda(t,y,k)\right] \,\pi(dk) dy\right|\, \le \, \Vert F\Vert_{L^2_\pi(\R\times \T)}\eps.
\]
Recalling \eqref{eq:prob_repr_bar_W} and \eqref{eq:prob_repr_tilde_W_la}, it is not difficult to check that Corollary
\ref{cor030203-20} now implies that
\[
\limsup_{\lambda\to+\infty}\int_{\R\times\T}F(y,k) \tilde{W}_\lambda(t,y,k)\, \pi(dk)dy\, = \, \int_{\R\times\T}F(y,k)
  \bar{W}(t-\delta,y)\, \pi(dk)dy.
\]
Choosing $\delta>0$ sufficiently small, we can then guarantee that
\[
\left|\int_{\R\times\T}F(y,k)\left[ \bar
  W(t-\delta,y)-\bar
  W(t,y)\right]\, \pi(dk)dy\right|\, < \, \eps.
\]
As a result, we conclude that for any $\eps>0$,
\[
\limsup_{\lambda\to+\infty}\left|\int_{\R\times\T}F(y,k)\left[
  W_\lambda(t,y,k)-\bar
  W(t,y)\right]\, \pi(dk)dy\right|\, \le  \,\eps,
\]
and   \eqref{conv1} immediately follows. The conclusion of  Theorem
\ref{main-thm} with respect to the Lebesgue measure on $\T$, as in
\eqref{conv}, immediately follows from \eqref{conv1}  and the fact that $\|F
\tilde{W}_\lambda(t )\|_\infty \le \|F
W_0\|_\infty $. The claim that $\bar
  W(t,y)$ is a weak solution of   \eqref{apr204} follows from
  Proposition \ref{prop:prob_repres_limit_sol}.\qed

\bigskip

\subsection*{Proof of Lemma \ref{lm010203-20}}
Fix an arbitrary $\eps>0$. We show first  that there exists $\delta>0$ such that
\begin{equation}
\label{011202-20}
\limsup_{\lambda\to+\infty}\| W_\lambda(\delta)-\hat{W}^o_\lambda(\delta)\|_{L^2_{\pi}(\R\times\T)}^2\, < \, \eps,
\end{equation}
where $
\hat{W}^o_\lambda(t,y,k) := \bbE\left[ W_0\left(y, K_\lambda^o( t,k)\right)\right]$. Indeed, let $\rho>0$ be sufficiently small so that
\begin{equation}
\label{022201-20}
\limsup_{\lambda\to+\infty}\int_{\{|y|<\rho\}\times\T}| W_\lambda(\delta,y,k)-
\hat{W}^o_\lambda(\delta,y,k)|^2\, \pi(dk)dy \,  < \, \eps.
\end{equation}
We then consider $\rho'\in(0,\rho)$ and conclude that for any $k$ in $\T_\ast$,
\begin{equation}
\label{010403-20}
\left\{\sup_{|y|\ge\rho}\sup_{t\in[0,\delta]}
|Y^o_\lambda( t,y,k)-y|\ge \rho'\right\}\,
= \, \left\{\sup_{t\in[0,\delta]}
|Y_\lambda( t,0,k)|\ge \rho'\right\},
\end{equation}
where $Y_\lambda(t,0,k)$ is the analogue of the process $Y_\lambda^o( t,0,k)$, but without the interface (or, equivalently, with $p_+\equiv 1$ in the model).
The above equation now implies that for any $\rho'< \rho$, there exists a sufficiently small $\delta>0$ such that
\begin{equation}
\label{042201-20}
\limsup_{\lambda\to+\infty}\,\bbP\left(\sup_{|y|\ge\rho}\sup_{t\in[0,\delta]}
|Y^o_\lambda( t,y,k)-y|\ge \rho'\right)\, < \, \eps\rho.
\end{equation}
In particular, it then follows that
\begin{equation}
\label{032201-20a}
\lim_{\delta\to0+}\limsup_{\lambda\to+\infty}\sup_{|y|\ge\rho}\bbP\left(\delta\ge
\mathfrak{s}^\lambda_{y,k,1}\right) \, = \, 0.
\end{equation}
where $\mathfrak{s}^\lambda_{y,k,1}$ is the first time the process $Y_\lambda(
  t,y,k)$ crosses the interface. From the above reasoning,  we now claim that there exist
$\rho,\delta>0$ sufficiently small so that
\begin{equation}
\label{052201-20}
\limsup_{\lambda\to+\infty}\int_{\{\rho\le |y|\le\rho^{-1}\}\times\T}\left[W_\lambda(\delta,y,k)-\hat{W}^o_\lambda(\delta,y,k)\right]^2\, \pi(dk)dy \, <
\, 2\Vert W_0 \Vert^2_\infty\eps.
\end{equation}
Indeed, the expression inside the integral in \eqref{052201-20} can be rewritten as
\begin{multline*}
\bigl\{\bbE\left[
 W_0\left(Y_\lambda^o( \delta,y,k), K_\lambda^o(
  \delta,k)\right)-W_0\left(y, K_\lambda^o(
  \delta,k)\right)\right]\bigr\}^2
\\
\le \, \bbE\left[\sup_{k'\in\T}\left(
  W_0\left(Y_\lambda^o( \delta,y,k), k'\right)-W_0\left(y, k'\right)\right)^2\right].
\end{multline*}
Using \eqref{042201-20}, we then choose $\delta>0$ so that for any
$k\in\T_*$ and any $|y|\ge \rho$, we have
\[
\limsup_{\lambda\to+\infty}\bbE\left[\sup_{k'\in\T}\left(
  W_0\left(Y_\lambda^o( \delta;y,k), k'\right)-W_0\left(y, k'\right)\right)^2\right] \, < \, \Vert W_0 \Vert^2_\infty\eps\rho.
\]
Noticing that the expression inside the integral in \eqref{052201-20}
is uniformly bounded in $\lambda$, we can finally use Fatou's lemma to
conclude that   \eqref{052201-20} holds. Since $W_0$ is compactly supported, Equation \eqref{011202-20}, now follows if we prove that for a fixed $\delta>0$, it holds that
\begin{equation}
\label{052201-20b}
\lim_{\rho\to0+}\limsup_{\lambda\to+\infty}\int_{\{ |y|\ge \rho^{-1}\}\times\T} W_\lambda^2(\delta,y,k)\, \pi(dk)dy \, = \, 0.
\end{equation}
To do so, we first notice that the function
\[
W_0^*(y)\, := \, \sup_{k\in\T}|W_0(y,k)|+\sup_{k\in\T}|W_0(-y,k)|
\]
is compactly supported, even, non-negative and such that $W_0^*\left(Y^o_\lambda(\delta,0,k)\right)\le W_0^*\left(Y_\lambda(\delta,0,k)\right)$. Noticing that a classical argument through stable central limit theorem (see, e.g.,  \cite[Theorem 4.1]{durrett-resnick}) shows that $Y_\lambda(\delta,0,k)$ weakly converges to $\eta(\delta)=\zeta(\bar{\theta}\delta)$, we get that
\[
 \begin{split}
\limsup_{\lambda\to+\infty}\int_{\{|y|\ge \rho^{-1}\}\times\T} W_\lambda^2(\delta,y,k)\, \pi(dk)dy \, &\le \, \limsup_{\lambda\to+\infty}
\bbE\left[\int_{\{|y|\ge \rho^{-1}\}}\left[W_0^*(y+Y_\lambda(
  \delta,0,k))\right]^2\, dy\right] \\
  &= \, \bbE\left[\int_{\{|y|\ge \rho^{-1}\}}\left[W_0^*(y+\eta(\delta))\right]^2 \, dy\right].
  \end{split}
\]
Equation \eqref{052201-20b} now follows immediately taking the limit, as $\rho$ goes to zero, in the above expression.
We can now use \eqref{022201-20}, \eqref{052201-20} and
\eqref{052201-20b} to conclude the proof of   \eqref{011202-20}. Using the fact that $W_0$ is compactly supported
and \eqref{032201-20a}, we know that
for any $\eps>0$ there exists $\delta>0$ such that
\begin{equation}
\label{021202-20}
\limsup_{\lambda\to+\infty}\int_{\R\times\T}\left[\hat{W}^o_\lambda(\delta,y,k)- \hat{W}_\lambda(\delta,y,k)\right]^2\, \pi(dk)dy\, < \, \eps,
\end{equation}
where
$ \hat{W}_\lambda(t,y,k):=  \bbE\left[W_0\left(y, K_\lambda( t,k)\right)\right]$.
Noticing that the dynamics of the momentum process $K_\lambda( t,k)$ is reversible
with respect to the measure $\pi$ on the torus $\T$ and $0$ is a simple
eigenvalue for the generator $L_k$, we finally have that
\begin{equation}
\label{021202-20aa}
\limsup_{\lambda\to+\infty}\int_{\R\times\T}\left[
\hat{W}_\lambda(\delta,y,k)-\bar W_0(y)\right]^2\, \pi(dk)dy \, <\, \eps
\end{equation}
and we have concluded the proof. \qed

\section{Proof of Theorem \ref{cor012402-20}}
\label{sec:Proof_converg}

To prove Theorem \ref{cor012402-20}, we show the convergence of
the Markov semigroups corresponding to processes
$\left\{Z_{\lambda}^o(t,y,k)\right\}_{t\ge0}$, see \eqref{eq:def_Z_o_lambda} . In the first part of
the present section, we focus on constructing a   L\'evy-type process
$\{\hat{Z}^o_\lambda(t,y)\}_{t\ge0}$ whose increments,  after the first jump,
coincide with those of $\left\{Z_{\lambda}^o(t,y,k)\right\}_{t\ge0}$.

\subsection{Construction of the associated Markov process} Let us denote by $r(y)$ the density
of $S(K_1)$, where $S(k)$ is defined in \eqref{SS} and $K_1$ is distributed according to $R_2(k)dk$ (see \eqref{barR}).
Recalling that the function $S(k)$ is bijective, the law of $S(K_1)$ is then given by
\begin{equation}
\label{eq:def_mu_S}
\mu_S(A)\, := \, \int_A r(y) \, dy\, = \, \int_{S(k)\in A}R_2(k)\, dk,\quad A\in {\cal B}(\R).
\end{equation}
Since $S$ is odd, its law $\mu_S$ is even, and thus $r(y)$ is even.
Let us now define, for simplicity,
\[
 \tilde p_{\iota,\lambda}(z)\, := \,
   p_{\iota}\left(S^{-1}(\lambda^{1/\alpha}z)\right), \quad\iota \in \{-1,0,1\}
\]
for any $\lambda>0$ and any $z$ in $\R$, with the classical notation
$\tilde p_{\iota}(z):=\tilde p_{\iota,1}(z)$. Clearly, all the above
functions are even. Thanks to  \eqref{g}  we   may assume, with no
loss of generality, that
\begin{equation}
  \tilde{p}_0\left(y\right)\succeq \frac{1}{(1+\log |y|)^{\kappa}},\quad
  |y|\ge 1\label{eq:11c1}.
\end{equation}
 From the general assumptions of the model (cf.
\eqref{g}, \eqref{tot}, \eqref{S0} and \eqref{S01}), it is not difficult to verify, by direct calculations, the following properties:
\begin{lemma}
\label{lm010309-19}
We have
\begin{equation}
\label{eq:11b}
  r\left(y\right)\, \approx \, \frac{1}{(1+|y|)^{ 1+\alpha }}, \qquad  y\in \R
\end{equation}
and
\begin{equation}
\label{022002-20}
\lim_{y\to+\infty}  r\left(y\right)y^{ 1+\alpha }\, = \, r_* \, := \, \frac{R^\ast_2S_\ast^{1+\alpha}}{S'_\ast}.
\end{equation}
In addition,
\begin{equation}
\label{eq:11c2}
\liminf_{y\to+\infty} \tilde{p}_0\left(y\right)\log^{\kappa} y=:\tilde p_*>0.
\end{equation}
where
$
\tilde p_*:= p_*\beta_3^{\kappa},
$ and $p_*$ is given by  \eqref{g}.
\end{lemma}

Consider, for any $y$ in $\R$,
\begin{equation}
\label{bar-r}
\bar{r}(y) \, := \, \int_0^{+\infty}r\left(\frac{y}{\tau}\right)\frac{e^{-\tau}d\tau}{\tau},
\end{equation}
and its rescaled version $\bar{r}_\lambda(y):=\lambda^{1+1/\alpha}
  \bar{r}\left(\lambda^{1/\alpha}y\right)$ for any $\lambda>0$. Lemma \ref{lm010309-19} then implies that
\begin{align}
\label{eq:11b-bar} \bar{r}(y) \, &\approx \, \frac{1}{(1+|y|)^{ 1+\alpha }} \quad \text{ on } \,\, \R; \\
\label{eq:11b-bar1}
\lim_{y\to+\infty} \bar{r}(y)y^{1+\alpha} \, &= \, \bar r_* \, = \, r_*\int_0^{+\infty}\tau^{\alpha}e^{-\tau}\, d\tau;\\
\label{eq:11bblm} \lim_{\lambda\to+\infty} \bar{r}_{\lambda}(y) \, &= \, \frac{\bar
  r_*}{|y|^{ 1+\alpha}}, \quad \mbox{if } y \neq 0.
\end{align}
We now consider the Markov process $\left\{\hat Z_\lambda^o(t,y)\right\}_{t\ge0}$ starting
at $y$, whose generator is defined on $B_b(\R)$ by ${\hat{ \mathcal{L}}}_{\lambda}^o u(0)=0$ and for any $y\in \R_\ast$, by
\begin{equation}
\label{eq:def_hat_Z_0_lambda}
{\hat{\mathcal{L}}}_{\lambda}^ou(y) \, =  \int_{\R}\hat r_{\lambda}(y,y')[u(y')-u(y)]\,dy' + k_{\lambda}(y)\left[u(0)-u(y)\right],
\end{equation}
where the jump and killing kernels are given by
\[\begin{split}
\hat r_{\lambda}(y,y') \, &:= \, \mathds{1}_{yy'>0}\bar r_{\lambda}(y'-y)+\mathds{1}_{yy'<0}\left[\tilde p_{+,\lambda}(y'-y)\bar r_{\lambda}(y'-y)+\tilde p_{-,\lambda}(y'+y)\bar r_{\lambda}(y'+y)\right];\\
k_{\lambda}(y)\, &:= \, \int_{y}^{+\infty}\tilde{p}_{0,\lambda}(z)\bar  r_{\lambda}(z)\, dz\, = \, \int_{|y|}^{+\infty}\tilde{p}_{0,\lambda}(z)\bar  r_{\lambda}(z)\, dz.
\end{split}
\]
We will prove at the end of the present section that the just constructed Markov processes $\hat{Z}^o_\lambda(t,y)$ converges as well to the limit process $\zeta^o(t,y)$ defined in \eqref{Zkr}.
\begin{thm}
\label{thm012102-20}
Let $y$ be in $\R_*$. As $\lambda\to+\infty$, the processes $\{\hat Z_{\lambda}^o(t,y)\}_{t\ge0}$ converge both in finite distributions and weakly, over $\mathcal{D}[0,+\infty)$ with the $J_1$-topology, to $\left\{\zeta^o(t,y)\right\}_{t\ge0}$.
\end{thm}

\subsection{Properties of the  Markov semigroup corresponding to $\{\hat{Z}^o_\lambda(t,y)\}_{t\ge0}$}

In order to show Theorem \ref{thm012102-20}, we will strongly rely on
a convergence property between the corresponding Markov
semigroups. The process $\hat Z_\lambda^o(t,y)$ killed at the
interface is Markovian and its transition semigroup is given by
\begin{equation}
\label{eq:def_hat_P_o_lambda}
P_t^{o,\lambda}u(y)\,:= \,\bbE\left[u\left(\hat Z_\lambda^o(t,y)\right),\,t<\hat{\mathfrak
    s}_{y,{\mathfrak f}}^{\lambda}\right], \qquad u\in B_b(\R).
\end{equation}
 Here $\hat{\mathfrak s}^{\lambda}_{y,{\mathfrak f}}\, := \, \inf\{t>0\colon \hat Z_\lambda^o(t,y) =0\}$.

We shall also consider the process stopped at $\hat{\mathfrak
  s}^{\lambda}_{y,{\mathfrak f}}$. Its transition semigroup equals
\begin{equation}
\label{eq:def_semigroup1}
\hat{P}_{t}^{o,\lambda}u(y)\, := \, \bbE\left[u\left(\hat Z_{\lambda}^o(t,y)\right)\right]\, = \, P_{t}^{o,\lambda}u(y)+u(0)\bbP\left(t\ge \hat{\mathfrak
    s}_{y,{\mathfrak f}}^{\lambda}\right).
\end{equation}
According to Theorem \ref{thm012102-20} the    process $\zeta^o(t,y)$
is the limit of the killed processes $\{\hat{Z}^o_\lambda(t,y)\}_{t\ge0}$.
Its   Markov semigroup  satisfies the following.
\begin{prop}
\label{prop012411-20}
For any $y$ in $\R_\ast$, the process $\{\zeta^o(t,y)\}_{t\ge0}$ generates a symmetric Markov semigroup $\{P_t^o\}_{t\ge0}$ on $L^2(\R)$ given by
\begin{equation}
\label{eq:def_P0_t}
P_t^ou(y)\, = \,\bbE\left[u(\zeta^o(t,y)),\,t<  \mathfrak{t}_{y,\mathfrak{f}}\right],\quad t\ge0.
\end{equation}
\end{prop}
\proof
Since the proof is quite long, we split it into two parts.

\subsubsection*{Proof that $(P^o_t)_{t\ge0}$ is a Markovian semigroup}
Fix $t,s>0$ and $N\in \N$. Let us consider  $0<s_1<\dots<s_N\le  s$, a
finite family $\{\phi_j, \,  j=1,\ldots,N \}$ of bounded Borel measurable functions $\phi_j\colon\R_*\to\R_+$ and
a bounded Borel measurable function $u\colon \R_*\to\R_+$. It is then enough to show that
\begin{equation}
\label{012411-20}
\bbE\left[u(\zeta^o(t+s,y))\Phi(y),\,t+s<  \mathfrak{t}_{y,\mathfrak{f}}\right]\, = \, \bbE\left[P_t^o u(\zeta^o(s,y))\Phi(y),\,s<  \mathfrak{t}_{y,\mathfrak{f}}\right].
\end{equation}
where we denoted, for simplicity, $\Phi(y)=\prod_{j=1}^N\phi_j(\zeta^o(s_j,y))$. Recalling the definition of $\zeta^o(t,y)$ in \eqref{Zkr}, we rewrite the left-hand side of \eqref{012411-20} as:
\begin{multline}
\label{012411-20a}
\sum_{0\le n\le m}\bbE\left[u(\zeta^o(t+s,y))\Phi(y),\, {\mathfrak t}_{y,n} \le s<  {\mathfrak t}_{y,n+1},\,
{\mathfrak t}_{y,m} \le t+s<  {\mathfrak t}_{y,m+1}\right]\\
=\, \sum_{0\le n\le
  m}\sum_{\eps_1,\ldots,\eps_{m}\in\{\pm1\}}\bbE\left[u\left(\zeta(t+s,y)\prod_{j=1}^{m}\eps_j\right)\Phi(y),\right. \\
 \sigma_1=\eps_1,\ldots,\sigma_m=\eps_m,\, {\mathfrak t}_{y,n} \le s<  {\mathfrak t}_{y,n+1},\,
{\mathfrak t}_{y,m} \le t+s<  {\mathfrak t}_{y,m+1}\Bigg].
\end{multline}
Let now $\{ \tilde\zeta(t)\}_{t\ge0}$ be  an independent copy of
the stable process $\zeta(t)$. Similarly to
\eqref{eq:def_stopping_times_frak_s}, we can also consider, for any
$m\in\{0,1,\ldots\}$  and  $z\in \R_\ast$, the $m$-th
consecutive time $\tilde{\mathfrak t}_{z,m}$
that the process $\{z+\tilde\zeta(t) \}_{t\ge0}$ crosses the interface.
Using the  independence of increments for the stable L\'evy process, we then rewrite
the right-hand side of  \eqref{012411-20a} as
\begin{multline}
\label{012511-20}
\sum_{0\le n\le
  m}\sum_{\eps_1,\ldots,\eps_{m}\in\{\pm1\}}\bbE\left[u\left(\left(\zeta(s,y)+\tilde\zeta(t)\right)\prod_{j=1}^{m}\eps_j\right)\Phi(y),\right.
  \\
\,\sigma_1=\eps_1,\ldots,\sigma_m=\eps_m,\, {\mathfrak t}_{y,n} \le s<  {\mathfrak t}_{y,n+1}  ,\,
\tilde{\mathfrak t}_{\zeta(s,y),m-n} \le t<  \tilde{\mathfrak t}_{\zeta(s,y),m-n+1}\Bigg].
\end{multline}
Then, using the symmetry of the law of a stable process, it follows that for any $\eps=\pm1$ and any $z$ in $\R_\ast$, we have that
\[\left\{\tilde{\zeta}(t),\{\tilde{t}_{z,j}\}_{j\ge0} \right\}_{t\ge
    0}\, \overset{\text{(law)}}{=} \,
  \left\{\eps\tilde{\zeta}(t),\{\tilde{t}_{\eps z,j}\}_{j\ge0}\right\}_{t\ge0}.\]
Therefore, we can finally rewrite \eqref{012511-20} as:
\[
\begin{split}
&\sum_{0\le n<
  m}\sum_{\eps_1,\ldots,\eps_{m}\in\{\pm1\}}\bbE\Bigg[u\left(\left(\zeta^o(s,y)+\tilde\zeta(t)\right)\prod_{j=n+1}^m\eps_j\right)\Phi(y),\\
&\qquad\qquad
\qquad \qquad\sigma_1=\eps_1,\ldots,\sigma_m=\eps_m,\, {\mathfrak t}_{y,n} \le s<  {\mathfrak t}_{y,n+1}  ,\,
\tilde{\mathfrak t}_{\zeta^o(s,y), m-n} \le t<  \tilde{\mathfrak
  t}_{\zeta^o(s,y), m-n+1}\Bigg]\\
&\,
=\, \sum_{n\ge 0}\sum_{\eps_1,\ldots,\eps_{n}\in\{\pm1\}}\bbE\left[ P_t^o u\big(\zeta^o(s,y)\big)\Phi(y),
\,\sigma_1=\eps_1,\ldots,\sigma_n=\eps_n,\, {\mathfrak t}_{y,n} \le s<  {\mathfrak t}_{y,n+1}\right]
\end{split}
\]
and   \eqref{012411-20} then follows immediately.

\emph{The semigroup $P^o_t$ is symmetric.} Fixed $u,v$ in $L^2(\R)$, our aim is to show the following:
\begin{equation}
\label{010112-20}
\int_{\R}P_t^ou(y)v(y)\, dy\, = \, \int_{\R}u(y)P_t^ov(y)\, dy, \qquad t>0.
\end{equation}
We start by rewriting the left-hand side of \eqref{010112-20} as
\begin{multline}
\label{020112-20}
\sum_{m\ge0} \sum_{\eps_1,\ldots,\eps_{n}\in\{\pm1\}}\int_{\R}\bbE\left[u\left(\zeta(t,y)\prod_{j=1}^m\eps_j\right)v(y),
\,\sigma_1=\eps_1,\ldots,\sigma_m=\eps_m,\,
{\mathfrak t}_{y,m} \le t<  {\mathfrak t}_{y,m+1}\right]\, dy\\
=\, \sum_{m\ge0}\sum_{\eps_1,\ldots,\eps_{m}\in\{\pm1\}}\prod_{j=1}^mp_{\eps_j}\int_{\R}\bbE\left[u\left(\zeta(t,y)\prod_{j=1}^m\eps_j\right)v(y),\,
{\mathfrak t}_{y,m} \le t<  {\mathfrak t}_{y,m+1}\right]\, dy.
\end{multline}
By a standard argument, using the Chapman-Kolmogorov equations, it is
possible to show the following
\begin{lemma}
\label{lemma:symmetry_semigroup}
Let $F\colon \mathcal{D}[0,t]\to \R_+$ bounded measurable. Then,
\begin{equation}
\label{010507-22}
\int_\R\bbE\left[u(\zeta(t,y))v(y)F(\zeta(\cdot,y))\right] \, dy\\
= \, \int_\R\bbE\left[u(y)v(\zeta
(t,y))F(\zeta^r(\cdot,y))\right] \, dy,
\end{equation}
where the reversed time process is given by
$\zeta^r(s,y):=\zeta((t-s)_{-},y)$, $s\in[0,t]$.
\end{lemma}
We can now exploit the above lemma with
$F(\zeta)=\mathds{1}_{\mathfrak{t}_{m}(\zeta) \le t<  \mathfrak{
    t}_{m+1}(\zeta)}$. Note that
$F(\zeta)=F(\zeta^r)$. Therefore,
we conclude  that
\begin{multline}
\label{eq:proof_symm}
\int_\R\bbE\left[u\left(\zeta(t,y)\prod_{j=1}^m\eps_j\right)v(y),\,
{\mathfrak t}_{m} \left(\zeta(\cdot,y)\right)\le t<  {\mathfrak t}_{m+1}\left(\zeta(\cdot,y)\right)\right] \, dy\\
= \, \int_\R\bbE\left[u\left( y\prod_{j=1}^m\eps_j\right)v(\zeta
(t,y)),\,
{\mathfrak t}_{m}\left(\zeta(\cdot,y)\right)  \le t< {\mathfrak t}_{m+1}\left(\zeta(\cdot,y)\right) \right] \, dy.
\end{multline}
If we assume for the moment that $\prod_{j=1}^m\eps_j=1$, it then obviously follows from \eqref{eq:proof_symm} that
\begin{multline}
\label{eq:proof_symm1}
\int_{\R}\bbE\left[u\left(\zeta(t,y)\prod_{j=1}^m\eps_j\right)v(y),\,
{\mathfrak t}_{y,m} \le t<  {\mathfrak t}_{y,m+1}\right] \, dy\\
= \, \int_{\R}\bbE\left[u(y)\, v\left( \zeta
(t,y)\prod_{j=1}^m\eps_j\right),\,
{\mathfrak t}_{y,m}  \le t< {\mathfrak t}_{y,m+1} \right]\, dy,
\end{multline}
where ${\mathfrak t}_{y,m}:={\mathfrak t}_{m}\left(\zeta(\cdot,y)\right)$.
Suppose now that $\prod_{j=1}^m\eps_j=-1$. Since the laws of $\{\zeta(t,y)\}_{t\ge0}$ and $\{-\zeta(t,-y)\}_{t\ge0}$ are identical, we can rewrite \eqref{eq:proof_symm} as
\begin{equation}
\begin{split}
\label{eq:proof_symm2}
\int_\R &\bbE\left[u\left(\zeta(t,y)\prod_{j=1}^m\eps_j\right)v(y),\,
{\mathfrak t}_{y,m} \le\, t<  {\mathfrak t}_{y,m+1}\right]\, dy\\
&\qquad \qquad \qquad\qquad \qquad=\, \int_\R\bbE\left[u(-y)v(-\zeta
(t,-y)),\,{\mathfrak t}_{-y,m}  \le t< {\mathfrak t}_{-y,m+1} \right] \, dy \\
&\qquad \qquad \qquad\qquad \qquad= \, \int_{\R}\bbE\left[u(y)v\left(\zeta
(t,y)\prod_{j=1}^m\eps_j\right),\,
{\mathfrak t}_{y,m}  \le t< {\mathfrak t}_{y,m+1} \right]\, dy.
\end{split}
\end{equation}
Applying \eqref{eq:proof_symm1}-\eqref{eq:proof_symm2} to \eqref{020112-20} above, we can finally conclude the proof. \qed

As in \eqref{eq:def_semigroup1}, it is now clear that for any $u$ in $B_b(\R)$, it holds that
\begin{equation}
\label{pto2}
\hat{P}_{t}^{o}u(y)\,  := \, \bbE\left[u\left(\zeta^o(t,y)\right)\right] = \,  P_{t}^{o}u(y)+u(0)\bbP\left(t\ge {\mathfrak t}_{y,{\mathfrak f}}\right),
\end{equation}
where $\hat{P}_{t}^{o}u(y)$ is the semigroup associated with
$\{\zeta^o(t,y)\}_{t\ge 0}$.
We will show in Section \ref{sec4} below that the following result holds:
\begin{thm}
\label{thm:conv_semigroups}
As $\lambda$ tends to $+\infty$, the semigroups $\{ P_{t}^{o,\lambda}\}_{t\ge0}$, defined in \eqref{eq:def_hat_P_o_lambda}, strongly converge in $L^2(\R)$, uniformly on compact intervals, to the semigroup
$\{P_{t}^o\}_{t\ge0}$ given in \eqref{eq:def_P0_t}.
\end{thm}
Thanks to the above result, we can now prove that the Markov
process $\hat Z_{\lambda}^o(t,y)$ converges to
$\zeta^o(t,y)$.

\subsection{Proof of Theorem \ref{thm012102-20}}

\subsubsection{Convergence of finite dimensional distributions}
Since the generalisation to finite distribution marginals is immediate, we will show only the convergence of the one-dimensional
distributions. Recalling the definition of the process $\{\hat{Z}^o_\lambda(t,y)\}_{t\ge 0}$ in \eqref{eq:def_hat_Z_0_lambda}, we start by considering the symmetric L\'evy process $\{\hat Z_{\lambda}(t)\}_{t\ge0}$ whose L\'evy symbol is given by $ \Psi_\lambda(\xi)$ where, thanks to \eqref{eq:11b},
\begin{equation}
\label{psin}
\Psi_\lambda(\xi)\, := \, 2\lambda^{1+1/\alpha}\int_{0}^{+\infty}{\sin}^2\left(\pi\xi y\right) \bar r\left(\lambda^{1/\alpha }y\right)\, dy
  \,\succeq  \, |\xi|^{\alpha}\int_{|\xi| \lambda^{-1/\alpha }}^{+\infty}\frac{\sin^2\left(\pi y\right)}{y^{ 1+\alpha }}\, dy.
\end{equation}
As before, we then denote $\hat{Z}_\lambda(t,y):=y+\hat{Z}_\lambda(t)$ for any $y\in \R^\ast$.
If we assume for the moment that $|\xi|
\ge \lambda^{1/\alpha}$, we have that
\begin{equation}
\label{theta-s}
\Psi_\lambda(\xi) \, \succeq \, \lambda(|\xi|\lambda^{-1/\alpha })^{\alpha}\int_{|\xi| \lambda^{-1/\alpha }}^{+\infty}\frac{\sin^2(\pi y)}{y^{ 1+\alpha }}\, dy \, \succeq \, \theta_*\lambda ,
\end{equation}
where
\begin{equation}
\label{eq:def_theta_ast}
\theta_*\, := \, \inf_{z\ge 1}\left\{ z^{\alpha}\int_z^{+\infty}\frac{\sin^2 (\pi y)}{y^{ 1+\alpha }}\, dy \right\} \, > \, 0.
\end{equation}
On the other hand, if $ |\xi| \le \lambda^{1/\alpha}$, we can write from \eqref{psin} that
\begin{equation}
\label{K}
\Psi_\lambda(\xi) \, \succeq \,  |\xi|^{\alpha}\int_{1}^{+\infty}\frac{\sin^2( \pi y)}{y^{ 1+\alpha }} \, dy \,  \succeq\,  \theta_\ast|\xi|^{\alpha}.
\end{equation}
Fixing $c>0$, we now notice from the above controls that $\Psi_\lambda(\xi)\succeq \theta_\ast>0$ for all $|\xi|\ge c$.
Then, it follows that the random variable $\hat{Z}_\lambda(t)$ admits a density $f_\lambda(t,\cdot)$ and
belongs to $
L^1(\R)\cap L^\infty(\R)$. Noticing that  the process  $\{\hat
  Z_{\lambda}^o(t,y)\}_{t\ge0}$ can be obtained from $\{\hat
  Z_{\lambda}(t,y)\}_{t\ge0}$ by performing a
  transformation of the trajectory (consisting in  transmission-reflection, or moving
  the phonon to the interface) described in Section
\ref{sec4.3}, we have in particular that
\begin{equation}
\label{eq:proof_conv_distr}
|\hat
  Z_{\lambda}^o(t,y)|\, \le \, |\hat
  Z_{\lambda}(t,y)|, \quad \mathbb{P}\text{-a.s.}
\end{equation}
Since we already know the weak convergence of the laws of $\Big(\hat
  Z_{\lambda}(t,y)\Big)$, it in turn implies the tightness of the laws of the random
variables $\hat Z_{\lambda}^o(t,y)$, as $\lambda\to+\infty$.
To conclude the proof, it is enough to show that
for any $y\in\R_*$ and any $\phi\in C^\infty_c(\R)$, we have
 \begin{equation}
\label{022102-20}
\lim_{\lambda\to+\infty}\bbE\left[\phi\Big(\hat  Z_{\lambda}^o(t,y)\Big)\right]\, = \,\bbE\left[\phi\Big(\zeta^o(t,y)\Big)\right].
\end{equation}
Let $\eps>0$ be arbitrary. Choosing $h\in(0,t]$ sufficiently small, we can guarantee that
\[
\bbP\left({\mathfrak t}_{y,1}\le h\right)+\limsup_{\lambda\to+\infty}\bbP\left(\hat{\mathfrak{s}}_{y,1}^{\lambda}\le h\right)\, < \, \eps,
\]
where ${\mathfrak t}_{y,1}$, $\hat{\mathfrak{s}}^\lambda_{y,1}$
represents the first time the processes crosses the interface
$o$ (cf.\ Equation
\eqref{eq:def_stopping_times_frak_s}). Recalling the definition of the
semigroups $\hat{P}^{o,\lambda}_t$, $\hat{P}^{o}_t$ in
\eqref{eq:def_semigroup1} and \eqref{eq:def_P0_t}, we then have that
\begin{align}
\label{022102-20a}
\limsup_{\lambda\to+\infty}\left|\bbE
\left[\phi\left(\hat
Z_{\lambda}^o(t,y)\right)\right]-
\bbE\left[\hat{P}^{o,\lambda}_{t-h}\phi\left(\hat Z_{\lambda}
(h,y)\right)\right]\right|\, &\preceq \, \eps \\
\label{022102-20z}
\left|\bbE\left[\phi\left(
   \zeta^o(t,y)\right)\right]-
  \bbE\left[\hat{P}^o_{t-h}\phi\left(
    \zeta(h,y)\right)\right]\right|\, &\preceq \, \eps.
\end{align}
By the above reasoning and \eqref{eq:def_semigroup1}, it now follows that
\begin{equation}
\label{072202-20}
\bbE\left[\hat{P}^{o,\lambda}_{t-h}\phi\left(\hat Z_{\lambda} (h,y)\right)\right] \, = \, \int_{\R} P^{o,\lambda}_{t-h}\phi(z+y)f_{\lambda}(h,z) \, dz+\phi(0) \int_{\R}\bbP\left(t\ge \hat{\mathfrak s}^{\lambda}_{y+z, {\mathfrak f}}\right) f_{\lambda} (h,z) \, dz.
\end{equation}
Thanks to Proposition \ref{thm:conv_semigroups} and \cite[Theorem 46.2]{gnedenko}, we know that $ P^{o,\lambda}_{t-h}\phi$ and $f_{\lambda}(h,\cdot)$
strongly converge in $L^2(\R)$, as $\lambda\to+\infty$,  to $ P^o_{t-h}\phi$ and $f(h,\cdot)$, the density
of $\zeta(h)$, respectively. Therefore,
\begin{equation}
\label{072202-20b}
\lim_{\lambda\to+\infty}\int_{\R}P^{o,\lambda}_{t-h}\phi(z+y)f_{\lambda}(h,t)\, dz \, = \, \int_{\R} P^o_{t-h}\phi(z+y)f(h,z) \, dz.
\end{equation}
To deal with the above limit we  use  the following lemma, whose proof is presented in
Section \ref{sec5.3.2}.
\begin{lemma}
\label{lm012102-20}
Let $f$ be in $L^1(\R)$. Then,
\begin{equation}
\label{032102-20}
\lim_{\lambda\to+\infty}\int_{\R}f(y)\bbP\left(t<\hat{\mathfrak
    s}^{\lambda}_{y,{\mathfrak f}}\right)\, dy \, = \, \int_{\R}f(y)\bbP\left(t<{\mathfrak t}_{y,{\mathfrak f}}\right)\, dy.
\end{equation}
\end{lemma}
Recalling that $f_{\lambda}(h,\cdot)$
converge to $f(h,\cdot)$ in $L^2(\R)$,  Lemma
\ref{lm012102-20} implies that
\begin{equation}
    \label{072202-20c}
    \lim_{\lambda\to+\infty}\int_{\R}\bbP\left(t\ge \hat{\mathfrak
s}^{\lambda}_{y+z,{\mathfrak f}} \right) f_{\lambda}(h,z)\, dz \, = \, \int_{\R}\bbP\left(t\ge {\mathfrak
t}_{y+z,{\mathfrak f}} \right) f(h,z) \, dz.
\end{equation}
Hence, from \eqref{072202-20}, \eqref{072202-20b}, \eqref{072202-20c}, and \eqref{pto2}, we conclude that
\[
\lim_{\lambda\to+\infty}\bbE\left[\hat{P}^{o,\lambda}_{t-h}\phi\left(\hat
    Z_{\lambda}(h,y)\right)\right]\, = \, \bbE\left[\hat{P}^o_{t-h}\phi\left(\zeta^o(h,y)\right)\right].
\]
Combining the above equation with \eqref{022102-20a} and \eqref{022102-20z}, we can easily conclude the proof of the convergence of finite dimensional distributions.

\subsubsection{Proof of Lemma \ref{lm012102-20}}
\label{sec5.3.2}
We are going to show \eqref{032102-20} only for $f$ in
$C^\infty_c(\R)$. The general statement of Lemma \ref{lm012102-20} then follows from a density argument. Fixed $M>0$, let us consider an even, smooth function
$\phi_M\colon \R\to [0,1]$ such that
\begin{equation}
\label{eq:test_function}
\phi_M(y) \, = \,
\begin{cases}
 1, &\mbox{ if } |y| \le M; \\
 0, &\mbox{ if } |y| \ge 2M.
 \end{cases}
\end{equation}
We can then define $\psi_M:=1-\phi_M$. Using \eqref{eq:proof_conv_distr}, we now show that for any $|y|\le M/2$, it holds that
\[
\bbE\left[\psi_M\left(\hat Z_{\lambda}^o(t,y)\right)\right] \, \le \, \bbP\left[|\hat Z^o_{\lambda}(t,y)|\ge M\right] \,   \le \,  \bbP\left[|\hat
    Z_{\lambda}(t)|\ge M/2\right].
\]
Clearly, a similar reasoning holds for $\zeta^o(t,y)$ as well.
Fixed $\eps>0$, we now choose $M>0$ and $\lambda_0>1$ large enough so that for any $\lambda\ge \lambda_0$, $|y|\le M/2$, it holds that
\[
\bbE\left[\psi_M\left(\hat{Z}^o_\lambda(t,y)\right)\right]+ \bbE\left[\psi_M\left(\zeta^o(t,y)\right)\right]\, \preceq\, \eps.
\]
Using again that $M$ is large enough so that $\text{supp} f\subseteq
[-M/2,M/2]$, we then conclude that
\[
\limsup_{\lambda\to+\infty}\left|\int_{\R}f(y)\left[\bbP\left(\hat{\mathfrak
    s}^{\lambda}_{y,{\mathfrak f}}>t\right)- P_{t,\lambda}^o\phi_M(y)\right] \, dy\right|+
\left|\int_{\R}f(y)\left[\bbP\left({\mathfrak t}_{y,{\mathfrak f}}>t\right)- P_{t}^o\phi_M(y)\right]\,dy\right|\, \preceq \, \eps.
\]
The above estimate and Proposition \ref{thm:conv_semigroups} finally imply that
\[
\limsup_{\lambda\to+\infty}\left|\int_{\R}f(y)\left[\bbP\left(\hat{\mathfrak
    s}^{\lambda}_{y,{\mathfrak f}}>t\right)-\bbP\left({\mathfrak t}_{y,{\mathfrak f}}>t\right)\right]\, dy\right| \, \preceq \, \eps,
\]
and Equation \eqref{032102-20}  immediately follows.
\qed

\subsubsection{Weak convergence in $J_1$-topology}

\label{sectight}
By the previous argument, it is enough to
prove tightness of the laws of $\{\hat Z_{\lambda}^o(t,y)\}_{t\ge0}$
over $\mathcal{D}[0,+\infty)$ with the $J_1$-topology. We shall make
use of the following notation. Given a function $f$ in $D[0,t_*]$, $t_\ast>0$ and $a<b$ in $[0,t_\ast]$  we denote $
\omega\left(a,b,f\right) := \sup\{|f(t)-f(s)|\colon a \le s<t<b\}$. We can then define the $D$-modulus of $f$ of step $\delta>0$ as:
\begin{equation}
\label{omp}
\omega'_{[0,t_*]}(\delta,f) \,:=\, \inf_{\mathcal{P}\in\mathds{I}_\delta}\max_{j=1,\dots,N}\omega\left(t_{j-1},t_j,f\right),
\end{equation}
where $\mathds{I}_\delta$ is composed of all the partitions $\mathcal{P}=\{t_0=t<\dots<t_N=t_\ast\}$ of $[0,t_\ast]$ such that $t_j-t_{j-1}\ge \delta$, for any  $j=1,\ldots,N$. It is not difficult to check that the laws of $\{\hat Z_\lambda(t,y)\}_{t\ge0}$ are tight over $\mathcal{D}[0,+\infty)$ with the $J_1$-topology. Theorem $13.2$ in
\cite{billingsley} then implies that for any $t_*>0$ and any $\eps>0$, it holds that
\begin{gather}
\label{012402-20}
\lim_{a\to+\infty}\sup_{\lambda\ge1}\bbP\left(\sup_{t\in[0,t_*]}|\hat
  Z_\lambda(t,y)|\ge a\right) \, = \, 0;\\
  \label{022402-20}
\lim_{\delta\to 0^+}\sup_{\lambda\ge1}\bbP\left(\omega'_{[0,t_*]}\left(\delta,\hat
  Z_\lambda(\cdot,y)\right)\ge \eps\right) \, = \, 0.
\end{gather}
Moreover, using  nested partitions of the time interval $[0,t_\ast]$,
one can easily show that
\[
\omega'_{[0,t_*]}\left(\delta, \hat Z_\lambda^o(\cdot,y)\right) \, \le \,  2\omega'_{[0,t_*]}\left(\delta,\hat Z_\lambda(\cdot,y)\right).
\]
Thanks to the above control and \eqref{eq:proof_conv_distr}, we
therefore conclude  that  \eqref{012402-20} and \eqref{022402-20} hold for $\hat Z_\lambda^o(t,y)$ as well and the tightness of the latter, as $\lambda\to+\infty$, immediately follows.
\qed

\subsection{The end of the proof of Theorem \ref{cor012402-20}}

The conclusion of Theorem \ref{cor012402-20} follows from  Theorem
\ref{thm012102-20} and Proposition \ref{prop012002-20}. Indeed,
suppose that $\{\hat Z_{\lambda}^o(t,y)\}_{t\ge0}$, $y\in \R_*$ is
the process defined in the previous section. Let $\tau_0$ be an
exp($1$) distributed random variable independent of the
process. Define also the random variable $\si$ that is independent of the
process and $\tau_0$ such that $\bbP[\si=\iota]=p_\iota(k)$, $\iota\in\{-1,0,1\}$.
Recall that $A^\lambda(y,k)$ is defined by \eqref{Ala} and $Z_1^\lambda(y,k)=y-\lambda^{-1/\alpha}S(k)\tau_0$.  Let
$A^\lambda_\iota(y,k):=\{yZ_1^\lambda(y,k)<0\text{ and }
\si^{\lambda}=\iota\} .$
Note that
\[
Z_{\lambda}^o(t,y,k) \, \stackrel{\rm law}{=} \,\left\{
\begin{array}{ll}
y& \mbox{ if }  t<\tau_0/\la,\\
0&\mbox{on }A^\lambda_0(y,k), \mbox{ if } \tau_0/\la\le t,\\
\hat Z_{\lambda}^o\Big(t-\tau_0/\la,-Z_1^\lambda(y,k)\Big)&\mbox{on }A^\lambda_{-1}(y,k) \mbox{ if } \tau_0/\la\le t,\\
\hat Z_{\lambda}^o\Big(t-\tau_0/\la, Z_1^\lambda(y,k)\Big)&\mbox{ if otherwise}.
\end{array}
\right.
\]
Convergence of finite dimensional distributions, claimed in Theorem
\ref{cor012402-20}, is then a consequence of Proposition
\ref{prop012002-20} and Theorem  \ref{thm:conv_semigroups}. Tightness
can be argued from tightness of $\Big(Z_{\lambda}(t,y,k)\Big)$
analogously to  Section \ref{sectight}.
\qed

\section{Proof of Theorem \ref{thm:conv_semigroups}}
\label{sec4}
We are going to show here the strong $L^2$-convergence of the semigroups $ P_{t}^{o,\lambda}$, defined in \eqref{eq:def_hat_P_o_lambda} and associated with the Markov processes $\hat{Z^o}_\lambda(t,y)$, to the semigroup
$\{P_{t}^o\}_{t\ge0}$, given in \eqref{eq:def_P0_t} and corresponding to
the process $\zeta^o(t,y)$. The main tool used in the proof is
the notion of the $\Gamma$-convergence of the Dirichlet forms
corresponding to the semigroups. Before going into the actual proof, we recall some basic facts on the subject.

\subsection{Basic notions on Sobolev spaces and Dirichlet forms}\label{s.Df}
Given $\beta$ in $(0,2)$, we define the following form:
\begin{equation}
\label{e.dE}
{\cal E}[u]\, := \, \frac{1}{2}\int_{\R^2}[u(y')-u(y)]^2q_\beta(y'-y)\, dydy'
\end{equation}
for any   Borel function $u\colon \R\to \R$.   We admit the possibility
that
the right-hand side of \eqref{e.dE} equals infinity.
Then,
$H^{\beta/2}(\R) :=  \left\{u\in L^2(\R)\colon {\cal E}[u] < \infty\right\}$ is a Hilbert space when endowed with the following norm
\begin{equation}
\label{eq:def_Sobolev_norm}
    \|u\|_{H^{\beta/2}(\R)}\, := \, \left( \|u\|_{L^2(\R)}^2+ {\cal E}[u]\right)^{1/2}.
\end{equation}
We have $H^{\alpha/2}(\R)\subset C_b(\R)$ when, as in our case,
$\alpha>1$ (see, e.g., \cite{fukushima}, equation $(1.4.33)$) and thus $u(0)$ is well defined. Moreover,   $C_c^\infty(\R)$ is dense in
$H^{\alpha/2}(\R)$ and it is not difficult to check that the form
${\cal E}$ is \textit{regular} on $H^{\alpha/2}(\R)$ in the sense of   \cite[ Example 1.4.1]{fukushima}.

Recalling that $\R_\ast:=\R\smallsetminus \{0\}$, we also consider the
family $H_0^{\beta/2}(\R_*)$ of functions given by the completion of
$C_c^\infty(\R_*)$ under the norm $\Vert \cdot
\Vert_{H^{\beta/2}(\R)}$. When $\alpha>1$, the space
$H_0^{\alpha/2}(\R_*)$ can be equivalently characterised as the set of
functions $u$ in $H^{\alpha/2}(\R)$ such that $u(0)=0$. For a proof of
this fact, see Corollary \ref{cor010309-19} below. It easily follows
that the form ${\cal E}$ is regular on $H_0^{\alpha/2}(\R_*)$. By \eqref{eq:interface_condition}, we also notice that the form $\hat{\cal E}$, given in \eqref{hat-C1}, is comparable with ${\cal E}$ on the space $H^{\alpha/2}_0(\R_*)$. In particular, it is regular and thus, it is a Dirichlet form on $H^{\alpha/2}_0(\R_*)$.

To prove Proposition \ref{thm:conv_semigroups}, we will show that the corresponding Dirichlet forms, defined below,
converge in a suitable sense. More precisely (cf. \cite[Definition
4.1]{DM93}, or  \cite[Section 1]{mosco}),
\begin{df}
\label{df2.5}
Let $\{\mathcal{E}_{\lambda}\}_{\lambda>0}$ be a family of Dirichlet forms on $L^2(\R)$ endowed with their natural domains:
\[
{\cal D}\left({\cal  E}_\lambda\right)\, :=\, \left\{u\in L^2(\R)\colon  {\cal  E}_\lambda[u]<+\infty\right\}.
\]
Then, $\mathcal{E}_\lambda$ is called $\Gamma$-convergent  to a Dirichlet form $\mathcal{E}_\infty$, as $\lambda\to+\infty$, if for any $u\in L^2(\R)$
the following conditions are satisfied:
\begin{itemize}
\label{item}
\item[i)]
for any family $\{u_{\lambda}\}_{\lambda>0}$ weakly convergent to $u$ in $L^2(\R)$,
it holds that
\[
\liminf_{\lambda\to+\infty}{\cal E}_{\lambda}[u_{\lambda}]\, \ge \,  {\cal E}_\infty[u]
\]
\item[ii)]
there exists a family $\{v_{\lambda}\}_{\lambda>0}$ strongly convergent to $u$  in $L^2(\R)$  such that
\[
\limsup_{\lambda\to+\infty}{\cal E}_{\lambda}[v_{\lambda}]\, \le \, {\cal E}_\infty[u].
\]
\end{itemize}
\end{df}
The notion of $\Gamma$-convergence is particularly useful for our purposes since it naturally implies (cf.\ \cite[Corollary $2.6.1$]{mosco}) the strong convergence of the corresponding semigroups on $L^2(\R)$.

We conclude this subsection presenting some properties of the Sobolev
spaces we have just introduced. The proofs of the results formulated
below can be found in Appendix \ref{secA}. Let us start with
the following variant of the Hardy-type inequality of Dyda \cite{dyda}.
\begin{prop}
\label{dyda}
For $\beta\not=1$, there exists a positive constant $C_\beta:=C(\beta)$ such that
\begin{equation}
\label{hardy}
\int_{\R}\frac{u^2(y)}{|y|^{\beta}}\, dy \, \le \,  C_\beta
\|u\|_{H^{\beta/2}_0 (\R_\ast)}^2,\quad \mbox{for all
}u\in H_0^{\beta/2}(\R_*).
\end{equation}
Moreover, if $1<\beta_0<\beta_1 <2$, then $\sup_{\beta\in[\beta_0,\beta_1]} C_\beta \, < \, +\infty$.
\end{prop}
We then obtain a characterisation of $H_0^{\beta/2}(\R_*)$ as a subspace of $H^{\beta/2}(\R)$.
\begin{prop}
\label{prop010209-19}
For $\beta\not=1$, we have
\begin{equation}
\label{020209-19}
H_0^{\beta/2}(\R_*)\, = \, \left\{u\in H^{\beta/2}(\R)\colon \int_{\R}\frac{u^2(y)}{|y|^{\beta}}dy<+\infty\right\}.
\end{equation}
\end{prop}
Finally, two other characterisations
of the space $H^{\alpha/2}_0(\R_\ast)$ are proposed for
indices $\alpha>1$.
\begin{cor}
\label{cor010309-19}
For every $\alpha$ in $(1,2)$, we have
\begin{align}
\label{010309-19}
H_0^{\alpha/2}(\R_*) \, = \, H^{\alpha/2}(\R)\cap
\bigcap_{1<  \beta <\alpha}H_0^{  \beta/2 }(\R_*); \\
\label{010309-19a}
H_0^{\alpha/2}(\R_*) \, = \, \left\{u\in H^{\alpha/2}(\R)\colon u(0)=0\right\}.
\end{align}
\end{cor}

\subsection{Properties of the associated Dirichlet forms}

Let us consider the Dirichlet form corresponding to $P_t^{o,\lambda}$:
\begin{equation}
\label{tC}
\hat{\cal E}_\lambda[u]\, := \, \frac{1}{2}\int_{\R^2}\hat
  r_{\lambda}(y,y')[u(y')-u(y)]^2\, dy'dy + \int_{\R} k_{\lambda}(y)u^2(y)\, dy,
\end{equation}
where the two kernels $\hat
  r_{\lambda}$, $k_{\lambda}$ were defined in \eqref{eq:def_hat_Z_0_lambda}. It follows from Proposition \ref{prop012411-20} that the process $\{\zeta^0(t,y)\}_{t\ge 0}$ admits a Dirichlet form $\mathcal{E}^{o}$ given by
\begin{equation}
\label{eq:cal_E_o}
\mathcal{E}^o[u]\, := \, \lim_{t\to0^+}\frac{1}{t}\int_{\R}\left[u(y)-P_t^ou(y)\right]u(y)\, dy,
\end{equation}
for any function $u$ belonging to its natural domain $\mathcal{D}(\mathcal{E}^o):=\{u\in L^2(\R)\colon \mathcal{E}^o[u]<\infty\}$. The main result of the present section is the following:
\begin{thm}
\label{prop012906-19}
The Dirichlet forms $\{\hat{\cal E}_{\lambda}\}_{\lambda>0}$  are $\Gamma$-convergent, as $\lambda$ goes to $+\infty$, to the Dirichlet form ${\cal E}^o$.
\end{thm}
 Before proving Theorem \ref{prop012906-19}, we show that $\mathcal{E}^o$ actually coincides  with the Dirichlet form $ \bar{r}_\ast\hat{\mathcal{E}}$, defined by \eqref{hat-C1}, and it is comparable to $\mathcal{E}$ with $\beta=\alpha$, which was given in \eqref{e.dE}.
\begin{prop}
\label{prop010212-20}
For any $u\in H_0^{\alpha/2}(\R_\ast)$, we have
    \begin{align}
    \label{equiv1}
    \mathcal{E}^o[u]\, &= \, \bar{r}_\ast\hat{\mathcal{E}}[u];\\
  \label{equiv}
\mathcal{E}^o[u] \, &\approx \, \mathcal{E}[u].
\end{align}
\end{prop}
\proof
For notational simplicity, we start by denoting
\begin{equation}
\label{eq:proof_Dirich_form}
\mathcal{E}^o_t[u]\, := \, \frac{1}{t}\int_{\R}[u(y)-
  P_t^ou(y)]u(y)\, dy,
\end{equation}
so that, in particular, $\mathcal{E}^o[u] = \lim_{t\to
  0^+}\mathcal{E}^o_t[u]$. Using the symmetry of the semigroup
$P^o_t$, we   write
\[
\begin{split}
\frac{1}{2t}\int_{\R} \bbE[(u(y)-u(\zeta^o(t,y)))^2 ,\,t<  \mathfrak{t}_{y,\mathfrak{f}}]\, dy\, &= \, \frac{1}{t}\int_{\R}P_t^ou^2(y)\, dy-\frac{1}{t}\int_{\R}u(y)P_t^ou(y)\, dy\\
&=\mathcal{E}^o_t[u]-\frac{1}{t}\left(\int_{\R}u^2(y)\big(1-P_t^o1(y)\big)\, dy\right).
\end{split}\]
We stress here that even though the constant function $1$ is not in
$L^2(\R)$, it is still possible to use the symmetry of the operator
$P_t^o$, arguing by approximation. Therefore,
\begin{equation}
\label{020412-20}
\mathcal{E}^o_t[u]\, = \, \frac{1}{t}\int_{\R}u^2(y)\mathbb{P}(t\ge \mathfrak{t}_{y,\mathfrak{f}})\, dy+\frac{1}{2t}\int_{\R} \bbE[(u(y)-u(\zeta^o(t,y)))^2 ,\,t<  \mathfrak{t}_{y,\mathfrak{f}}]\, dy.
\end{equation}
It is then easy to verify that, for any $u\in C_c^\infty(\R_*)$,
\begin{equation}
\label{eq:proof_Dirich}
\mathcal{E}^o[u] \, = \, \lim_{t\to 0^+}\mathcal{E}^o_t[u] \, = \,  \bar{r}_\ast\hat{\mathcal{E}}[u].
\end{equation}
Fixed $y$ in $\R_\ast$, let $\{\tilde\zeta(t)\}_{t\ge0}$ be the
symmetric $\alpha$-stable L\'evy process starting at $y$ but killed at
hitting $0$. Clearly, its transition semigroup $\{\tilde
P_t\}_{t\ge0}$, given as in \eqref{eq:def_P0_t}, is made of symmetric Markov contractions
on $L^2(\R)$. In particular, its corresponding Dirichlet form equals:
\[
\tilde{\cal E}[u] \, := \, \frac{1}{2}\int_{\R}\int_{\R}[u(y')-u(y)]^2
q_\alpha(y'-y)\, dydy' \, = \, \mathcal{E}[u],
\]
for any function $u$ belonging to its natural domain
$\mathcal{D}(\tilde{\mathcal{E}})=H_0^{\alpha/2}(\R_*)$. For a proof
see, e.g., \cite[Section 3.3.3]{MR2849840}.
We consider as well
$\tilde{\cal E}_t[u]$ defined as in \eqref{eq:proof_Dirich_form} with
respect to $\tilde{P}_t$. The same calculations that lead to
\eqref{020412-20} can be performed again to show that
\[
\tilde{\mathcal{E}}_t[u]\, = \, \frac{1}{t}\int_{\R}u^2(y)\mathbb{P}(t\ge \mathfrak{t}_{y,\mathfrak{f}})\, dy+\frac{1}{2t}\int_{\R} \bbE[(u(y)-u(\tilde{\zeta}(t,y)))^2 ,\,t<  \mathfrak{t}_{y,\mathfrak{f}}]\, dy.
 \]
Recalling the construction of the process $\zeta^o(t,y)$ in \eqref{Zkr}, we now write that
\begin{align}
\label{010412-20}
&
t\mathcal{
  E}^o_t[u] \, = \,\int_{\R}u^2(y)\mathbb{P}(t\ge \mathfrak{t}_{y,\mathfrak{f}})\, dy\\
&+\frac{1}{2}\sum_{m\ge0}\sum_{\eps_1,\ldots,\eps_{m}\in\{\pm 1\}}\prod_{j=1}^m p_{\eps_j}\int_{\R} \bbE\left[\left(u(y)-u\left(\zeta(t,y)\prod_{i=1}^m\eps_i\right)\right)^2 ,\,
{\mathfrak t}_{y,m} \le t<  {\mathfrak t}_{y,m+1}\right]\, dy.\notag
\end{align}
We then denote by ${\cal P}_m$ the family of sets
$\{\eps_1,\ldots,\eps_{m}\}\subseteq\{-1,1\}^m$ such that $\prod_{j=1}^m\eps_j=1$
and by ${\cal P}_m^c$ its complement. Let
\[
s_m\, := \, \sum_{\{\eps_1,\ldots,\eps_{m}\}\in {\cal P}_m}\prod_{j=1}^mp_{\eps_j}.
\]
In particular, $s+1=p_+$.
By a direct calculation, it is possible to show the following:
\begin{lemma}
\label{lm010412-20}
The sequence $\{s_m\}_{m\in \N}$ tends to $1/2$. Moreover, if $p_+>1/2$, then $\{s_m\}_{m\in \N}$ is strictly decreasing, if $p_+=1/2$, then it holds that $s_m=1/2$ for
any $m\ge1$ and if $p_+<1/2$, then $\{s_{2m-1}\}_{m\in \N}$ increases while $\{s_{2m}\}_{m\in \N}$ decreases.
\end{lemma}
From   \eqref{010412-20}, we now have that
\begin{align}
\notag
t\mathcal{
  E}^o_t[u] \, &\ge \, \int_{\R}u^2(y)\mathbb{P}(t\ge \mathfrak{t}_{y,\mathfrak{f}})\, dy+\frac{1}{2}\sum_{m\ge0}s_m\int_{\R} \bbE\left[\left(u(y)-u(\zeta(t,y)\right)^2 ,\,
{\mathfrak t}_{y,m} \le t<  {\mathfrak t}_{y,m+1}\right] \, dy\\ \label{010412-20a}
& \ge \, \left(p_+\wedge \frac{1}{2}\right)t\tilde {\cal
  E}_t(f).
\end{align}
If we suppose now that $u$ belongs to ${\cal D}(\mathcal{E}^o)$, the
closure of $C^\infty_c(\R_\ast)$ with respect to the form
$\mathcal{E}^o$, then \eqref{010412-20a} implies that $u$ belongs to
$\mathcal{D}(\tilde{\mathcal{E}})=H_0^{\alpha/2}(\R_*)$, as well. On the other hand, if $u$ is in $C_c^\infty(\R_*)$ then we can  use the Hardy inequality \eqref{hardy} in \eqref{eq:proof_Dirich} to show that $\mathcal{E}^o[u] \preceq \tilde{\mathcal
  {E}}[u]$. Such an estimate then extends to $H_0^{\alpha/2}(\R_*)$ and implies in
particular that $H_0^{\alpha/2}(\R_*)\subseteq \mathcal{D}(\mathcal{E}^o)$. We have thus proven that $H_0^{\alpha/2}(\R_*)=  \mathcal{D}(\mathcal{E}^o)$ and there exists a positive constant $C$ such that
\[
C^{-1}\tilde{\mathcal{E}}[u] \,  \le \,  \mathcal{E}^o[u] \, \le \,  C
\tilde{\mathcal{E}}[u],\quad u\in H_0^{\alpha/2}(\R_\ast),
\]
which ends the proof of the proposition.  \qed

\subsection{Proof of Theorem \ref{prop012906-19}}
In the present section, we are going to
show that the Dirichlet forms
$\{\hat {\cal E}_{\lambda}\}_{\lambda>0}$, defined in
\eqref{tC}, are $\Gamma$-convergent, when
$\lambda$ tends to $+\infty$, to the
Dirichlet form $\bar{r}_\ast \hat {\cal E}$ given in
\eqref{hat-C1}-\eqref{eq:11b-bar1}. Thanks to Proposition \ref{prop010212-20}, Theorem \ref{prop012906-19} will then follow immediately.

Recalling the meaning of the $\Gamma$-convergence in Definition \ref{df2.5}, condition ii) easily follows choosing the trivial family $\{v_\lambda\}_{\lambda> 0}$ of functions given by $v_\lambda=u$ and using the Lebesgue dominated convergence theorem, together with \eqref{p-plus1} and \eqref{eq:11bblm}. We now focus on showing condition i).  Let $\{u_{\lambda}\}_{\lambda>0}$ be a family of functions weakly convergent to $u$ in $L^2(\R)$. We start by noticing that if $\liminf_{\lambda\to+\infty}\hat {\cal E}_{\lambda}[u_{\lambda}]=+\infty$,
then condition i) clearly holds. We can then suppose without loss of generality that
\begin{equation}
\label{022708-19}
\liminf_{\lambda\to+\infty}\hat{\cal E}_{\lambda}[u_{\lambda}] \, < \, +\infty.
\end{equation}
Let $\{\lambda_n\}_{n\in \N}$ be a sequence in $(0,+\infty)$ such that $\lambda_n\to +\infty$ and
\[
\lim_{n\to+\infty}\hat {\cal E}_{n}[u_{n}]\, = \, \liminf_{\lambda\to+\infty}\hat {\cal E}_{\lambda}[u_{\lambda}],
\]
where we denoted $u_n:=u_{\lambda_n}$, $\hat {\cal E}_{n}:=\hat {\cal
  E}_{\lambda_n}$. We will use the following lemma, whose proof is
presented in   Section \ref{sec6.4}.
\begin{lemma}
\label{cor020309-19a}
Let $\{u_n\}_{n\ge1}$ be a bounded sequence in $L^2(\R)$ such that $\lim_{n\to \infty} \hat{\mathcal{E}}_n[u_n]<\infty$. Then, there exists a subsequence $\{u_{n_k}\}_{k\ge1}$ that is a.s.\ convergent to $u$.
\end{lemma}
We   claim that $u$ is in $H_0^{\alpha/2}(\R_*)$. Indeed, by Lemma
\ref{cor020309-19a} above, there exists a sub-sequence of
$\{u_n\}_{n\ge1}$, that for simplicity we  denote again by the same symbol, that is a.s.\ convergent in $\R$.
From \eqref{eq:def_hat_Z_0_lambda} and \eqref{p-plus1}, we now notice that
\[\hat{r}_{\lambda_n}(y,y') \, \ge \,   \mathds{1}_{yy'>0}\bar r_{\lambda_n}(y'-y)+\mathds{1}_{yy'<0}\tilde p_{+,\lambda_n}(y'-y)\bar r_{\lambda_n}(y'-y)\, \,  \succeq \,  \bar{r}_{\lambda_n}(y'-y).\]
Fatou's lemma and   \eqref{eq:11bblm} then imply that
\[
\lim_{n\to+\infty}\hat {\cal E}_n[u_n]\, \ge \, \frac{1}{2}
\int_{\R^2}\liminf_{n\to+\infty}\left[\hat{r}_{\lambda_n}(y,y')
  \left(u_n(y')-u_n(y)\right)^2\right] \, dy'dy\, \succeq \,  \mathcal {E}[u].
\]
From \eqref{022708-19}, we now have that $u\in H^{\alpha/2}(\R)$. Using again \eqref{022708-19}, we notice that
\begin{equation}
\label{eq:control_M_conv}
\int_{\R}u_n^2(y) k_{\lambda_n}(y)\, dy \, \preceq \, 1, \quad n\ge N,
\end{equation}
for some $N\in\N$. Assuming that  $\lambda^{1/\alpha }_n|y|\ge 1$, we can use \eqref{eq:11c1}-\eqref{eq:11b-bar} to show that
\begin{equation}
\label{011407-22}
 k_{\lambda_n}(y) \, \succeq \, \lambda_n\int_{\lambda_n^{1/\alpha}|y|}^{+\infty}
  \frac{\left[1+\log
 (1+ z)\right]^{-\kappa}}{\left(1+z\right)^{1+\alpha}}\, dz\,
\succeq \,\frac{\left[1+\log
  (1+\lambda_n^{1/\alpha }|y|)\right]^{-\kappa}}{ \left(\lambda_n^{-1/\alpha }+|y|\right)^{\alpha}}.
\end{equation}
The second case, where $\lambda^{1/\alpha}_n|y|< 1$, is trivial, since then
\[
 k_{\lambda_n}(y) \, \succeq \, \lambda_n\int_{1}^{+\infty}
  \frac{\left[1+\log
 (1+ z)\right]^{-\kappa}}{\left(1+z\right)^{1+\alpha}}\, dz\,\succeq \la_n.
\]
Thus, $k_{\lambda_n}(y)$  can also be controlled from below by the term
appearing on the utmost right hand side of \eqref{011407-22}. From
\eqref{eq:control_M_conv} and \eqref{011407-22}, it   follows that
 \begin{equation}
 \label{011309-19}
 \int_{\R}\frac{u_n^2(y)}{
  \left(\lambda_n^{-1/\alpha }+|y|\right)^{\alpha}}\cdot\frac{dy}{\left[1+\log
  (1+\lambda_n^{1/\alpha } |y|)\right]^{\kappa}} \, \, \preceq \,  1,
 \end{equation}
for any $n\ge N$.
We shall use   the following.
\begin{lemma}
\label{lm011309-19}
Let $\{u_n\}_{n\ge1}$ be a bounded sequence in $L^2(\R)$ such that $\lim_{n\to \infty} \hat{\mathcal{E}}_n[u_n]<\infty$ and the inequality \eqref{011309-19} holds. Then, there exists $\rho'$ in $(0,1)$ such that
\begin{equation}
 \label{021309-19}
  \sup_{n \ge1}\int_{\R}\frac{u_n^2(y)  dy}{
  (\lambda_n^{-1/\alpha }+|y|)^{ \alpha (1-\rho)}} \,\, < \, +\infty, \qquad \rho\in(0,\rho').
  \end{equation}
\end{lemma}
The lemma is shown in Section \ref{sec6.5}. We proceed with its
application to the proof of Theorem \ref{prop012906-19}.
Using Fatou's lemma in \eqref{021309-19}, we conclude that
\[
\int_{\R}\frac{u^2(y)}{|y|^{\alpha(1-\rho)}}\, dy \, < \, +\infty,
\]
for any $\rho\in(0,1)$. Thanks to Proposition \ref{prop010209-19} and Corollary  \ref{cor010309-19}, we then infer that $u$ is in
$H^{ \alpha /2}_0(\R_*)$. Since $\{u_n\}_{n\ge1}$ is a.s.\ convergent
to a function $u$ in $H^{ \alpha /2}_0(\R_*)$, we can use Fatou's lemma and \eqref{eq:11bblm} to write that:
\[
\liminf_{n\to+\infty}\hat {\cal E}_{n}[u_n]\, \ge \,
\frac{1}{2}\int_{\R^2}\lim_{n\to+\infty}(u_n(y')-u_n(y))^2\hat
  r_{\lambda_n}(y,y')\, dydy'\,
= \,  \bar{r}_\ast \hat  {\cal E}[u],
\]
and we have proven part i) of Definition \ref{df2.5}, which ends the
proof of Theorem \ref{prop012906-19}. %\qed

 What remains to be done is to prove the auxiliary results presented above.
Before doing so, we however need the following result:
\begin{lemma}
\label{lm012708-19}
Let $\{u_n\}_{n\ge 1}$ be bounded in $L^2(\R)$ such that $\lim_{n\to \infty} \hat{\mathcal{E}}_n[u_n]<\infty$. Then,
\[
\lim_{K\to+\infty}\int_{|\xi|\ge K}|\hat u_n(\xi)|^2 \, d\xi \, = \, 0,
\]
uniformly in $n$. Here we have denoted by $\hat f(\cdot)$ the Fourier
transform of a function $f\in L^2(\R)$.
\end{lemma}
\proof
Thanks to \eqref{p-plus1}, we know that $p_\ast := \inf_{k\in\T} p_+(k) \, > \, 0$.
Since $\lim_{n\to \infty} \hat{\mathcal{E}}_n[u_n]<\infty$, we then
conclude from   \eqref{tC} that
\[
1\, \succeq \, \hat {\cal E}_n[u_n] \, \ge \,
\frac{p_*}{2}\int_{\R^2}\left[u_n(y')-u_ n(y)\right]^2 \bar{r}_{\lambda_n}(y'-y) \, dydy',
\]
where $\bar r_{\lambda}$ is given by \eqref{bar-r}. The right-hand side is then of the same order of magnitude as
\[
\int_{\R}|\hat u_n(\xi)|^2d\xi\Big(\int_{\R}\sin^2\left(\pi\xi
  y\right) \bar  r_{\lambda_n}\left(y\right)\, dy\Big) \, \preceq \, 1.
\]
Fixing $K>0$, it now follows that
\begin{equation}
\label{092708-19}
\int_{|\xi|\ge K}\Psi_n(\xi)|\hat u_{\lambda_n}(\xi)|^2 \, d\xi\, \preceq \,  1,
\end{equation}
where $\Psi_n:=\Psi_{\lambda_n}$ was defined in \eqref{psin}. Thus,
for any $n$ so large that $ K\le \lambda_n^{1/\alpha}$, the estimates \eqref{092708-19}, \eqref{theta-s} and \eqref{K} imply that
\[
K^{\alpha} \int_{\lambda_n^{1/\alpha }\ge|\xi|\ge K}|\hat
  u_{\lambda_n}(\xi)|^2\, d\xi\, \preceq \, K^{\alpha} \int_{\lambda_n^{1/\alpha }\ge|\xi|\ge K}|\hat
  u_{\lambda_n}(\xi)|^2\, d\xi + \lambda_n \theta_*\int_{\lambda_n^{1/\alpha }\le |\xi|}|\hat
  u_{\lambda_n}(\xi)|^2\, d\xi \,  \preceq \, 1,
\]
and the conclusion of  Lemma \ref{lm012708-19}   follows. \qed

\subsection{Proof of Lemma \ref{cor020309-19a}.}
\label{sec6.4}
Fix $\delta>0$ and let us consider   the function $\phi_\delta$ given in \eqref{eq:test_function}.
We  claim that the sequence $\{u_n\phi_\delta\}_{n\ge1}$ is strongly
compact in $L^2(\R)$. Using Pego Criterion, see \cite[Theorem 3]{pego}, it is enough to show that
\begin{align}
\label{042708-19a}
\lim_{K\to+\infty}\int_{|y|>K}|u_n(y) \phi_\delta(y)|^2 \, dy\, = \, 0;
\\
\label{042708-19}
\lim_{K\to+\infty}\int_{|\xi|>K}|\widehat{u_n\phi_\delta}(\xi)|^2\, d\xi\, = \, 0,
\end{align}
uniformly in $n \in \N $. Clearly, \eqref{042708-19a} holds since the sequence is uniformly bounded in $L^2(\R)$. To show \eqref{042708-19}, we notice that
\begin{equation}
\label{052708-19}
\int_{|\xi|>K}|\widehat{u_n\phi_\delta}(\xi)|^2\, d\xi \, = \, \int_{|\xi|>K}\left|\int_{\R}\widehat{u}_n(\xi-\eta)\widehat{\phi}_\delta(\eta)\, d\eta\right|^2\, d\xi \,
\le \,  I_1(K,L)+I_2(K,L),
\end{equation}
where for any fixed $L>0$, we have
\[
\begin{split}
 I_1(K,L)\, &:= \, 2 \int_{|\xi|>K}\left|\int_{|\eta|\le
 L}\widehat{u}_n(\xi-\eta)\widehat{\phi}_\delta(\eta)\, d\eta\right|^2\, d\xi,\\
I_2(K,L)\, &:= \, 2 \int_{|\xi|>K}\left|\int_{|\eta|> L}\widehat{u}_n(\xi-\eta)\widehat{\phi}_\delta(\eta)\, d\eta\right|^2\, d\xi.
\end{split}
\]
The Cauchy-Schwartz inequality and Lemma \ref{lm012708-19} then imply that for any fixed $L>0$,
\begin{equation}
\label{062708-19}
\lim_{K\to+\infty}\sup_{n\in \N} I_1(K,L)\, \le \, 4L\|\phi_\delta\|_{L^2(\R)}^2\lim_{K\to+\infty}\sup_{n\in \N}\int_{|\xi|>K-L}\left|\widehat{u}_n(\xi)\right|^2\,d\xi \,  = \, 0.
\end{equation}
On the other hand, recalling that $\{u_n\}_{n\ge1}$ is uniformly
bounded in $L^2(\R)$ and $\hat \phi_\delta$ belongs to the Schwartz
class, we conclude that for any $\eps>0$ it is possible to choose $L:=L(\eps)$, sufficiently large, such that
\begin{equation}
\label{072708-19}
I_2(K,L)\, \preceq \,  2\Vert u_n\Vert^2_{L^2(\R)}\left(\int_{|\eta|>L}|\widehat{\phi}_\delta(\eta)|d\eta\right)^2 \, \preceq \, \eps,
\end{equation}
uniformly in $n$. From \eqref{052708-19}-\eqref{072708-19}, we can
then conclude that   \eqref{042708-19} follows. From the compactness of $\left\{u_n\phi_\delta\right\}_{n\ge1}$, we can now choose a subsequence that is a.e.\
convergent on $[-\delta,\delta]$. Using the Cantor diagonal argument,  we then
find a subsequence of  $\left\{u_n\right\}_{n\ge1}$ that converges a.e.\ on
$\R$. \qed

\subsection{Proof of Lemma \ref{lm011309-19}.}
\label{sec6.5}
First observe that since the sequence $\{u_n\}_{n\ge1}$ is bounded in $L^2(\R)$, it is enough to show that for any $\rho$ small enough, we have:
\begin{equation}
\label{021309-19a}
 \limsup_{n \to +\infty}\int_{-1}^1\frac{u_n^2(y)}{
  \left(\lambda_n^{-1/\alpha }+|y|\right)^{ \alpha (1-\rho)}}\, dy \,<+\infty.
\end{equation}
We will actually prove an analogue of \eqref{021309-19a} with the
integral over $[0,1]$, as the argument in the  case of $[-1,0]$ is similar.
The proof is divided into three separate steps.

\emph{Step 1.} We start by claiming that
\begin{equation}
 \label{021309-19b}
 \int_{0}^{c_n(\gamma)}\frac{u_n^2(y)\, dy}{
  (\lambda_n^{-1/\alpha}+y)^{\alpha (1-\rho)}}\, \preceq \, c_n(\gamma/2),
 \end{equation}
 for any $\lambda_n\ge e$, any $\rho$ in $(0,1)$ and any $\gamma>1$, where
\begin{equation}
\label{eq:def_c_n_gamma}
c_n(\gamma) \, := \, \exp\left\{-\log^{\gamma}_2\lambda_n\right\}
\end{equation}
and $\log_2x:=\log\log x$, $x>1$.
Indeed, let us denote
\[
I_m \, := \, \left\{y\in[0,1]\colon e^{m-1} \lambda_n^{-1/\alpha}\le
y<e^m\lambda_n^{-1/\alpha}\right\},
\]
for any $m\in \N$ and $I_{0}:=\left\{y\in[0,1] \colon
y<\lambda_n^{-1/\alpha}\right\}$. Estimate \eqref{011309-19} now implies that
\begin{equation}
\label{031309-19n}
  \sum_{m=0}^{N_1}a_m\int_{I_m}\frac{u_n^2(y)}{\left(\lambda_n^{-1/\alpha }+y \right)^{ \alpha (1-\rho)}}\, \preceq \,\sum_{m=0}^{N_1}\frac{1}{(1+m)^{\kappa}}\int_{I_m}\frac{u_n^2(y)}{\left(\lambda_n^{-1/\alpha }+y\right)^{\alpha}} \, dy \, \preceq \, 1,
 \end{equation}
where a positive integer $N_1 $ is such that
 $e^{N_1-1}\lambda^{-1/\alpha }_n \le 1 < e^{N_1}\lambda^{-1/\alpha}_n$ and
\[
a_m\, := \, \frac{e^{\rho\alpha(N_1-m) }}{(1+m)^{\kappa}}.
\]
It is not difficult to check now that the sequence
$\{ a_m\}_{m\ge 1}$
is decreasing and thus, its minimum in $1\le m\le
m_*(\gamma):=N_1-[\log^{\gamma}N_1]$ equals $a_{m_*(\gamma)}$. Hence,
\eqref{031309-19n} implies that
\[
a_{m_*(\gamma)}\int_{0}^{e^{m_*(\gamma)-N_1}}\frac{u_n^2(y)}{\left(\lambda_n^{-1/\alpha
    }+y\right)^{ \alpha (1-\rho)}}\, dy \,  \preceq \,  1,\quad n\in \N
\]
since
$e^{m_\ast(\gamma)}\lambda_n^{-1/\alpha}=e^{m_\ast(\gamma)-N_1}$. Estimate
\eqref{021309-19b} follows from an observation that $e^{m_*(\gamma)-N_1} \succeq c_n(\gamma)$ and
\[
\frac{1}{a_{m_*(\gamma)}} \, \preceq \, \frac{N_1^k}{e^{\rho \alpha[\log^\gamma N_1]}} \, \preceq \, \exp\left\{\kappa\log N_1-\rho\alpha\log^{\gamma}N_1\right\}\, \preceq \,  c_n(\gamma/2).\]
\emph{Step 2.} We then show that for any $\lambda_n>1$, the function $u_n$ can be decomposed as:
\begin{equation}
\label{051309-19}
u_n(y ) \, = \, u^{(1)}_n(y)+u_n^{(2)}(y),
\end{equation}
for two functions $u^{(1)}_n$, $u^{(2)}_n$ such that
\begin{align}
\label{051309-19a}
\sup_{y\in\R}|u^{(1)}_n(y+h)-u^{(1)}_n(y)| \, &\preceq \,
  h^{(\alpha-1)/2},\quad  h>0,\\
  \label{051309-19b}
  \|u^{(2)}_n\|_{L^2(\R)}\, &\preceq \, 1/\sqrt{\lambda_n}.
\end{align}
We start by defining
\[
u^{(1)}_n(y)\, := \, \int_{|\xi|\le
  \lambda^{1/\alpha }}\exp\left\{2\pi i \xi y\right\}\hat
  u_n(\xi)\, d\xi,\qquad
u^{(2)}_n(y)\, := \, \int_{|\xi|>
  \lambda^{1/\alpha }}\exp\left\{2\pi i \xi y\right\}\hat
  u_n(\xi)\, d\xi.
\]
Clearly,   \eqref{051309-19} then holds.
We now use again the function $\Psi_n$ defined in \eqref{psin}. Fixed
$h>0$, we apply  \eqref{K} to write that
\[
|u^{(1)}_n(y+h)-u^{(1)}_n(y)|\,
\preceq \, \frac{1}{\sqrt{\theta_\ast}}\int_{|\xi|\le
\lambda_n^{1/\alpha }}\Psi_n^{1/2}(\xi)\left|\hat u_{\lambda}(\xi)\right|\frac{\left|1-\exp\left\{2\pi
  i \xi h\right\}\right|}{|\xi|^{ \alpha /2}}\, d\xi.
\]
Using the Cauchy-Schwartz inequality and then \eqref{092708-19}, we
conclude that estimate \eqref{051309-19a} holds:
\[
|u^{(1)}_n(y+h)-u^{(1)}_n(y)| \,
 \preceq \, \left(\int_{\R}\frac{\left|1-\exp\left\{2\pi
  i \xi \right\}\right|^2}{|\xi|
 ^{ \alpha }}\, d\xi\right)^{1/2}h^{(\alpha-1)/2}\, \preceq \,  h^{(\alpha-1)/2}.\]
On the other hand, estimate  \eqref{051309-19b} follows from
\eqref{theta-s} and \eqref{092708-19}:
\[
\|u^{(2)}_n\|_{L^2(\R)}^2\,
\le \, \frac{1}{\theta_\ast\lambda_n}\int_{|\xi|>
\lambda_n^{1/\alpha }}\Psi_n(\xi)\left|\hat u_{\lambda}\left(\xi\right)\right|^2\, d\xi \,
\preceq \,  \frac{1 }{\lambda_n}.
\]
\emph{Step 3.}   Fix $n$ large enough so that $\lambda_n\ge e$ and $\gamma>1$, we consider
\[
\tilde I_m \, := \, \{y\in[0,1]\colon mc_n(\gamma) \le
y<(m+1)c_n(\gamma)\},
\]
for any $m\in \N_0$. We can then write
\begin{equation}
 \label{071309-19}
\int_{0}^1\frac{u_n^2(y)}{
  (\lambda_n^{-1/\alpha}+y)^{ \alpha (1-\rho)}} \, dy \, \le \,  \sum_{m=0}^{N_2} \int_{\tilde I_m}\frac{u_n^2(y)}{
  (\lambda_n^{-1/\alpha }+y)^{ \alpha(1-\rho)}} \, dy \, \preceq \,
  I_n+I\!I_n
 \end{equation}
where $N_2$ is the (unique) positive integer such that $N_2 c_n(\gamma)\le
1< (N_2+1)c_n(\gamma)$ and
\[
\begin{split}
I_n \, := \, \sum_{m=0}^{N_2} \int_{0}^{c_n(\gamma)}\frac{u_n^2(y)}{(\lambda_n^{-1/\alpha }+mc_n(\gamma)+y)^{\alpha (1-\rho)}}\, dy;\\
I\!I_n \,:= \, \sum_{m=1}^{N_2} \int_{0}^{c_n(\gamma)}\frac{[u_n(y+mc_n(\gamma))-u_n(y)]^2}{
  (\lambda_n^{-1/\alpha }+mc_n(\gamma)+y)^{\alpha(1-\rho)}}\, dy.
\end{split}
\]
To control the first term $I_n$, we start by noticing that if $|y|\le c_n(\gamma)$, then $
\lambda_n^{-1/\alpha }+y \le 2c_n (\gamma)$ for sufficiently large $n $. We can then write:
\[
\begin{split}
\sum_{m=0}^{N_2} \frac{1}{(\lambda_n^{-1/\alpha }+mc_n(\gamma)+y)^{\alpha (1-\rho)}}\,
\preceq \,
  \frac{1}{(\lambda_n^{-1/\alpha }+y)^{\alpha (1-\rho)}}+\frac{1}{c_n^{ \alpha (1-\rho)}
  (\gamma)} \,
\preceq \, \frac{1}{(\lambda_n^{-1/\alpha }+y)^{ \alpha (1-\rho)}},
\end{split}\]
provided that $\rho$ is so small that $ \alpha (1-\rho)>1$.
From \eqref{021309-19b} we thus have that
$I_n\preceq c_n(\gamma/2)\preceq 1$, uniformly in $n$. Concerning the second term $I\!I_n$, we use \eqref{051309-19} to decompose it as $I\!I^{(1)}_n+I\!I^{(2)}_n$ where for any $j=1,2$, $I\!I^{(j)}_n$ is obtained from $I\!I_n$ by replacing there $u_n$ with $u^{(j)}_n$. Noticing that $
N_2= \left[c^{-1}_n(\gamma)\right]$,
Estimate \eqref{051309-19a} now implies that
\begin{align}
\label{011409-19}
I\!I^{(1)}_n \,\preceq \,
 \sum_{m=1}^{N_2}
  \frac{c_n(\gamma)(mc_n(\gamma))^{\alpha-1}}{(mc_n(\gamma))^{\alpha (1-\rho)}} \, \preceq \,
  c_n^{\alpha\rho}(\gamma)\sum_{m=1}^{N_2}\frac{1}{m^{1-\rho \alpha}} \, \preceq \, (c_n(\gamma)
  N_2)^{ \alpha \rho} \, \preceq \, 1,
\end{align}
 To control $I\!I^{(2)}_n$, we recall  that
$\lambda_n   $ dominates
$c_n^{-\alpha (1-\rho)}(\gamma) =  \exp\left\{ \alpha
  (1-\rho)\log_2^{\gamma}\lambda_n \right\}$ (cf. \eqref{eq:def_c_n_gamma})
for $n$ sufficiently large,
and then,
thanks to \eqref{051309-19b}, we write
\[
\begin{split}
I\!I^{(2)}_n \, &\preceq \,
 \Vert u^{(2)}_n\Vert^2_{L^2(\R)}\sum_{m=1}^{N_2}
 \frac{1}{(mc_n(\gamma))^{ \alpha (1-\rho)}} \, \preceq \, \frac{1}{\lambda_n
  c_n^{\alpha (1-\rho)}(\gamma)}\, \preceq \, 1.
\end{split}
\]
This ends the proof of  \eqref{021309-19a}.  \qed

\appendix
\section{Proofs of the results of Section \ref{s.Df}}
\label{secA}

\subsection{Proof of Proposition \ref{dyda}}

If we assume   that $u\in C_c^\infty(\R_*)$, then inequality
\eqref{hardy} follows from \cite[Theorem $1.1$ part $T4)$]{dyda}. By Fatou's lemma, we can then extend \eqref{hardy} to all functions $u$ in $H_0^{\beta/2}(\R_*)$ with the same constant.
The constant $C_\alpha$ is bounded on compact subsets of $(1,2)$, which follows from the  inspection of the constant in the last inequality of the
proof of \cite[Lemma $3.3$]{dyda}, see also  the proof of \cite[Theorem 1.1]{dyda}.
\qed

\subsection{Proof of Proposition \ref{prop010209-19}} For simplicity, let us denote the set on the right hand side of
\eqref{020209-19} by ${\cal H}$. We start by noting that Hardy inequality \eqref{hardy} immediately implies that $
H_0^{\beta/2}(\R_*)\subset {\cal H}$. In order to show the other
inclusion, let us take $u\in \mathcal{H}$. For any $\delta>0$, let us
consider again the function $\phi_\delta$, as considered in
\eqref{eq:test_function}, and define
$\psi_\delta=1-\phi_\delta$. Without any loss of generality, we can also assume that $\|\psi_\delta'\|_{\infty}\le 2/\delta$. Let $u_\delta:=u\psi_\delta$. It is easy to check that $u_\delta$ is in $\mathcal{H}$ and that it can be approximated by elements of $C_c^\infty(\R_*)$ in $H^{\beta/2}(\R)$-norm.
To conclude, it is then enough to show that $\lim_{\delta\to 0}\Vert v_\delta \Vert_{H^{\beta/2}(\R)} = 0$ where we denoted $v_\delta:=u_\delta-u$. Clearly, $\Vert v_\delta \Vert_{L^2}\to 0$ as $\delta$ goes to $0$, by the dominated convergence theorem. On the other hand, $\mathcal{E}[v_\delta]$ can be decomposed as
\[ c_\beta\iint_{yy'<0}\frac{[v_{\delta}(y')-v_{\delta}(y)]^2}{|y'-y|^{1+\beta}}\, dydy' + c_\beta \iint_{yy'>0}\frac{[v_{\delta}(y')-v_{\delta}(y)]^2}{|y'-y|^{1+\beta}}\, dydy' \, =: \, I_\delta+I\!I_\delta.
\]
To control the first term $I_\delta$, we notice that for any fixed $\varepsilon>0$, there exists $\delta:=\delta(\epsilon)>0$ sufficiently small such that
\begin{equation}
\label{proof:eq:control_I_delta}
I_\delta\, \preceq \,
\iint_{yy'<0}\frac{[v^2_{\delta}(y')+v^2_{\delta}(y)]}{|y'-y|^{1+\beta}}dydy'\, \preceq \, \iint_{y
  y'<0}\frac{v^2_{\delta}(y)}{|y'-y|^{1+\beta}}\, dydy'\, \preceq  \, \int_{-2\delta
  }^{2\delta}\frac{u^2(y)}{|y|^{\beta}} \, dy\, < \, \varepsilon.\end{equation}
To deal with the second term $I\!I_\delta$, we start by splitting it as
\[
\iint_{yy'>0,\, |y-y'|\ge 4\delta}\frac{[v_{\delta}(y')-v_{\delta}(y)]^2}{|y'-y|^{1+\beta}}\, dydy'+
\iint_{yy'>0,\, |y-y'|< 4\delta} \frac{[v_{\delta}(y')-v_{\delta}(y)]^2}{|y'-y|^{1+\beta}}\, dydy' \, =: \, I\!I_\delta^1+I\!I_\delta^2.\]
Recalling that $v_\delta(y) = 0$ if $|y|\ge 2\delta$, we denote $A=\{yy'>0, |y-y'|\ge 4\delta, |y|\le
  2\delta\}$ so that
\begin{equation}
\label{proof:eq:control_II_delta_1}
I\!I_\delta^1\, \preceq \,
\iint_A\frac{v^2_{\delta}(y)}{|y'-y|^{1+\beta}}\, dydy'\, \preceq \, \int_{ |y|\le
  2\delta}\frac{v^2_{\delta}(y)}{(4\delta)^{\beta}}\, dy \,
\preceq\, \int_{ |y|\le
  2\delta}\frac{v^2_{\delta}(y)}{|y|^{\beta}}\, dy \, < \, \varepsilon,
\end{equation}
provided that $\delta:=\delta(\epsilon)>0$ is sufficiently small. Similarly, we also have that
\[
I\!I_\delta^{2}\,
\le \, \iint_{yy'>0,\, |y|\vee|y'|< 10\delta}\frac{[v_{\delta}(y')-v_{\delta}(y)]^2}{|y'-y|^{1+\beta}}\, dydy \, \preceq \,  I\!I_\delta^{2,1}+I\!I_\delta^{2,2}+I\!I_\delta^{2,3},
\]
where we have denoted
\[
\begin{split}
I\!I_\delta^{2,1}\, &:= \, \iint_{|y|\vee|y'|< 6\delta} \frac{[\psi_{\delta}(y')-\psi_{\delta}(y)]^2u^2(y')}{|y'-y|^{1+\beta}}\, dydy'\\
 I\!I_\delta^{2,2}\, &:=\, \iint_{|y|\vee|y'|< 6\delta}
\frac{[u(y')-u(y)]^2 \psi_{\delta}^2(y)}{|y'-y|^{1+\beta}}\, dydy',\\
I\!I_\delta^{2,3}\, &:= \, \iint_{|y|\vee|y'|< 6\delta} \frac{[u(y')-u(y)]^2}{|y'-y|^{1+\beta}}\, dydy'.
\end{split}
\]
Recalling that   $u\in \mathcal{H}$, it follows that for any sufficiently small $\delta>0$,
\begin{equation}
\label{proof:eq:control_II_delta_22}
I\!I_\delta^{2,2}+I\!I_\delta^{2,3}\,< \,\varepsilon.
\end{equation}
Since we have chosen $\psi_\delta$ such that $\|\psi_\delta'\|_{\infty}\le 2/\delta$, we have that
\begin{equation}
\label{proof:eq:control_II_delta_13}
I\!I_\delta^{2,1}\, \preceq \,
\delta^{-\beta}\int_{ |y'|< 6\delta} u^2(y')\,  dy'\,
\preceq \,  \int_{ |y'|< 6\delta
  }\frac{u^2(y')}{|y'|^{\beta}} \, dy'\, <\,\varepsilon
\end{equation}
for sufficiently small $\delta>0$. Summarising,  estimates
\eqref{proof:eq:control_I_delta} --
\eqref{proof:eq:control_II_delta_13} imply that $\mathcal{E}[v_\delta]<\epsilon$ for any $\delta$ small enough, and the conclusion of the proof then follows.
\qed

\subsection{Proof of Corollary \ref{cor010309-19}}
First, we show  \eqref{010309-19}. For simplicity, let us denote the
set appearing on its  right-hand side by $\mathcal{H}'$. Then, the
inclusion $H_0^{\alpha/2}(\R_*)\subset \mathcal{H}'$ trivially follows
once we   show that
\begin{equation}
\label{011509-19}
\|u\|^2_{H^{\beta/2}(\R)}\le\left(1+
\frac{8}{\beta}
\right)\|u\|^2_{H^{\alpha/2}(\R)},
\end{equation}
for any $\beta<\alpha$ in $(0,2)$ and $u$ in $H^{\alpha/2}(\R)$. If
$u\in C^\infty_c(\R)$, then the above inequality immediately follows from
\[
\iint_{|y-y'|\ge1}\frac{[u(y')-u(y)]^2}{|y'-y|^{1+\beta}}\,dydy' \, \le\,  2\iint_{|y-y'|\ge1}\frac{u^2(y')+u^2(y)}{|y'-y|^{1+\beta}}\, dydy'\,\le\,  \frac{8}{\beta}\|u\|_{L^2(\R)}^2.
\]
The general conclusion can be then obtained by an approximation
argument. To show the other inclusion, let us take now a function $u$
in  $\mathcal{H}'$. Proposition \ref{dyda} and estimate
\eqref{011509-19} then imply that for any $(1+\alpha)/2\le \beta\le
\alpha$, we have
\[
\int_{\R}\frac{u^2(y)}{|y|^{\beta}}\, dy \, \le \, \left(\sup_{(1+\alpha)/2\le \beta\le \alpha}C_\beta\right) \left(1+\frac{8}{\beta}\right)\|u\|_{H^{\alpha/2}(\R)}^2.
\]
Thanks to Fatou's lemma and Proposition \ref{prop010209-19}, it
follows that $u\in H_0^{\alpha/2}(\R_*)$. We have thus showed   \eqref{010309-19}.
To prove   \eqref{010309-19a}, let us suppose that $u\in
H^{\alpha/2}(\R)$ satisfies  $u(0)=0$. From \eqref{010309-19} and Proposition \ref{prop010209-19}, it suffices to show that
\[
\int_{\R}\frac{u^2(y)}{|y|^{\beta}}\, dy\, < \, +\infty,
\]
for any $\beta\in (1,\alpha)$. Since $\alpha>1$, we recall from
\cite[Theorem 1, p. 394]{muramatu} that $u$ possesses an H\"older continuous representative
and that in particular,
\[
\sup_{y\neq y'\in\R}\frac{|u(y)-u(y')|}{|y-y'|^{(\alpha-1)/2}} \, \le \, C \|u\|_{H^{\alpha/2}(\R)}.
\]
It then follows that
\[
\int_{\R}\frac{u^2(y)}{|y|^{\beta}}\, dy \, = \, \int_{\R}\frac{(u(y)-u(0))^2}{|y|^{\beta}}\, dy \,
\preceq\, \int_{|y|\le1}\frac{|y|^{\alpha-1}}{|y|^{\beta}}\, dy+\int_{|y|\ge1}u^2(y)\, dy \,
< \, +\infty,
\]
and the proof of the corollary is  concluded.
\qed

\section{Proof of Proposition \ref{prop:prob_repres_limit_sol}}
\label{appC}

\subsection{Representation of a weak solution.}
We will assume without loss of generality that $T_\infty=0$. We also
suppose  for the moment that $T_o=0$. Then Proposition
\ref{prop:prob_repres_limit_sol} can be obtained form the observation that
\[\bar{W}(t,y)\, = \, \bbE\left[\bar{W}_0(\eta^o(t,y))\right] \, = \, P^{o}_{\bar{\theta}t}\bar{W}_0(y),\]
is the semigroup solution to Cauchy problem \eqref{apr204}, where
$P^o_t$ and $\bar{\theta}$ have been defined in \eqref{eq:def_P0_t} and \eqref{bar-theta}, respectively.

If $T_o\neq 0$, let us consider again the smooth function $\phi_M$ given in \eqref{eq:test_function}. Since $\bar{W}_0-T_o\phi_M$ is in $L^2(\R)\cap \mathcal{H}_o$ and $\big(\bar{W}_0-T_o\phi_M\big)(0)=0$, we already know that
\[\bar{W}_M(t,y)\,  := \, \bbE\left[\left(\bar{W}_0-T_o\phi_M\right)(\eta^o(t,y))\right] \, = \, P^{o}_{\bar{\theta}t}\left(\bar{W}_0-T_o\phi_M\right)(y),\]
is a weak solution of \eqref{apr204} for $T_o=0$. Thus, $\bar{W}_M\in C\left([0,+\infty);L^2(\R)\right)\cap L^2_{\rm
  loc}\left([0,+\infty);  \mathcal{H}_o\right)$ and
\begin{multline}
\label{021803-19xy}
\int_{\R}F(0,y)[ \bar{W}_0(y)-T_o\phi_M(y)]\, dy\,
  =\, \int_{\R}F(t,y)[\bar{W}_M(t,y)]\, dy\\
-\int_0^{t}\int_{\R}\partial_sF(s,y) [\bar{W}_M(s,y)]\, dy ds +\bar{\gamma} \int_0^{t}\hat{\cal
  E} [F(s,\cdot),\bar{W}_M(s,\cdot)] \, ds.
 \end{multline}
It is not difficult to check now that $\bar{W}_0-T_o\phi_M\to
\bar{W}_0-T_o$, as $M\to+\infty$, both in $\mathcal{H}_o$ and pointwise. Recalling that $P^o_t$ is the Markov semigroup associated with the form $\hat{\mathcal{E}}$, we have that
\[\hat{\mathcal{E}}\left[ \bar{W}_M(t)\right]\, = \, \hat{\mathcal{E}}\left[P^{o}_{\bar{\theta}t}\left(\bar{W}_0-T_o\phi_M\right)\right] \, \le \, \hat{\mathcal{E}}\left[\bar{W}_0-T_o\phi_M\right].\]
It then follows that $\bar{W}_M\to\bar{W}-T_o$ in $\mathcal{H}_o$ and $\bar{W}-T_o$ is in $L^2_{\text{loc}}(0,+\infty,\mathcal{H}_o)$. Finally, we can follow the same arguments in the proof of Lemma \ref{lm012102-20} to show that
\[\lim_{M\to \infty}\int_\R F(y)\bar{W}_M(t,y) \, dy \, = \, \int_\R F(y)\left[\bar{W}(t,y)-T_o\right] \, dy,\]
for any $F$ in $C^\infty_c(\R)$. Hence, we can pass to the limit in
\eqref{021803-19xy} and we conclude that $\bar{W}(t,y)$ given by
\eqref{eq:prob_repr_bar_W} is a weak solution in the sense of
Definition \ref{df011803-19}.

\subsection{Uniqueness of a weak solution.}
For notational simplicity, let us  denote by $\langle\cdot,\cdot\rangle_{\mathcal{H}_o}$ the scalar product on $\mathcal{H}_o$. Let $\bar{W}_1$, $\bar{W}_2$ be two weak solutions (in the sense of Definition \ref{df011803-19}) of Cauchy problem \eqref{apr204}. By definition, we know that for each $j=1,2$, there exists $T^j_\infty$ such that $\bar{W}_j-T^j_\infty\in
C([0,+\infty);L^2(\R))$ and $\bar{W}_1(0)=\bar{W}_2(0)=\bar{W}_0$. It then follows immediately that $T^1_\infty=T^2_\infty$. This in turn implies that
\[
\bar{W} \, := \, \bar W_1-\bar{W}_2 \, \in \,  C\left([0,+\infty);L^2(\R)\right)\cap L^2_{\rm
  loc}\left([0,+\infty);  \mathcal{H}_o\right).
\]
Assuming for the moment that $\bar{W}$ belongs to $C^1\left([0,+\infty);  \mathcal{H}_o\right)$, it is easy to check that
\begin{equation}
\label{022101-22}
\left\|\int_0^{t}\bar W (s)\,ds\right\|^2_{\mathcal{H}_0}\, = \, 2\int_0^t
  \left\langle \bar W (s),\int_0^{s}\bar W (s')\, ds'\right\rangle_{\mathcal{H}_o} \, ds, \quad t\ge0.
\end{equation}
In the more general case where $\bar{W}$ is in $L^2_{\rm
  loc}\left([0,+\infty); \mathcal{H}_o \right)$, \eqref{022101-22} follows by a mollification argument in time. Recalling that $\bar{W}$ satisfies \eqref{021803-19xx} for $\bar{W}_0=T_o=0$, we can now deduce that for any $F \in  L^2(\R)\cap\mathcal{H}_o$, it holds that
\begin{equation}
\label{012101-22}
\int_0^{t}\langle F,\bar
  W(s)\rangle_{\mathcal{H}_o}\, ds\, = \,-\langle F,\bar W(t)\rangle_{L^2(\R)}.
 \end{equation}
In particular, from \eqref{022101-22} and \eqref{012101-22} with $F=\bar{W}(s)$, we finally get that
\[
\label{042101-22}
\left\|\int_0^{t}\bar W (s)\,ds\right\|^2_{\mathcal{H}_o}\, = \, -2\int_0^t\|\bar W (s)\|_{L^2(\R)}^2\, ds,
\]
which implies that $\bar W (t)=0$.
\qed

\end{document}